# Computational multiphase micro-periporomechanics for dynamic shear banding and fracturing of unsaturated porous media

Hossein Pashazad, Xiaoyu Song[*]

*Engineering School of Sustainable Infrastructure and Environment*
*University of Florida, Gainesville, FL, USA*

Abstract

Dynamic shearing banding and fracturing in unsaturated porous media is a significant problem in engineering and science. This article proposes a multiphase micro-periporomechanics (µPPM) paradigm for modeling dynamic shear banding and fracturing in unsaturated porous media. Periporomechanics (PPM) is a nonlocal reformulation of classical poromechanics to model continuous and discontinuous deformation/fracture and fluid flow in porous media through a single framework. In PPM, a multiphase porous material is postulated as a collection of a finite number of mixed material points. The length scale in PPM that dictates the nonlocal interaction between material points is a mathematical object that lacks a direct physical meaning. As a novelty, in the coupled µPPM, a microstructure-based material length scale is incorporated by considering micro-rotations of the solid skeleton following the Cosserat continuum theory for solids. As a new contribution, we reformulate the second-order work for detecting material instability and the energy-based crack criterion and J-integral for modeling fracturing in the µPPM paradigm. The stabilized Cosserat PPM correspondence principle that mitigates the multiphase zero-energy mode instability is augmented to include unsaturated fluid flow. We have numerically implemented the novel µPPM paradigm through a dual-way fractional-step algorithm in time and a hybrid Lagrangian-Eulerian meshfree method in space. Numerical examples are presented to demonstrate the robustness and efficacy of the proposed µPPM paradigm for modeling shear banding and fracturing in unsaturated porous media.

## 1. Introduction

Dynamic shearing banding and fracturing in unsaturated deformable porous media (e.g., geomaterials and human tissues) is a significant problem in engineering and science (e.g., [1–7]). The coupled shear banding/cracking and fluid flow processes can significantly deteriorate the integrity of civil infrastructure and damage soft/hard human tissues under certain circumstances (e.g., [4–6, 8–12]). For instance, fault propagation triggered by earthquakes involves coupled shear banding/cracking and fluid flow in geomaterials (e.g., [13–16]). Landslides and embankment failure could be triggered by shear banding or cracking in geomaterials [17–21]. Computational coupled poromechanics is an essential tool in studying the coupled localization or fracturing and fluid flow processes in porous materials (e.g., [1, 2, 7, 12, 22]). Coupled periporomechanics (PPM) [23–32] is a strong nonlocal reformulation of classical poromechanics [33, 34]. In PPM, the coupled motion equation and mass balance equation are formulated in terms of integral-differential equations (integration in space and differentiation in time) instead of partial differential equations [29] through the peridynamic state and effective force state concepts [27, 35, 36]. For comparison between PPM and other continuum-based numerical methods, we refer to the literature [29]. Through the stabilized multiphase correspondence principle, classical constitutive models and physics laws can be incorporated in PPM for modeling the coupled deformation, shear banding/fracturing, and fluid flow processes in porous media [30, 37]. In PPM, a multiphase porous medium is postulated as a collection of a finite number of mixed material points. The length scale in PPM that dictates the nonlocal interaction between material points is a mathematical object that lacks a clear physical meaning. In the present study, as a new contribution, we formulate a robust multiphase micro-periporomechanics (μPPM) for modeling localized failure and fracturing in unsaturated porous media in which the length scale is related to micro-rotations of material points. In this new coupled μPPM, different from the previously formulated PPM, the length scale is associated with the micro-structure of porous media in line with the classical micro-polar continuum theory [38, 39] and thus the non-locality has a rational physical meaning. It is noted that the micro-polar/Cosserat continuum mechanics [38, 39] was applied to model shear banding in solids and porous media by the computational geomechanics community several decades ago. The pioneering work on the subject can be found in the distinguished literature (e.g., [40–45]). For instance, the micro-polar continuum theory was adopted to regularize the finite element solution of strain localization in granular materials [41, 42, 46]. The micro-polar poromechanics has been developed to model strain localization in saturated porous media [43, 44]. Meanwhile, the micro-polar continuum theory has been used to model fracture in solids [45]. For instance, in [45], the authors formulated a micro-polar J-integral [47, 48] for modeling cracking. Recently, the micro-polar continuum theory [38, 39] has been used in peridynamics to model fracturing in solids [49, 50]. For stance, a state-based elastic PD model was proposed in [50] to model cracks in brittle materials. For a review of the micro-polar PD for solids, we refer to [51]. In [51], the authors formulated a viscous-plastic Cosserat PPM for modeling dynamic shear bands and crack branching in single-phase porous media in which the micro-polar length scale was incorporated following the Cosserat continuum theory for solids [38, 39]. This study extends the single-phase Cosserat PPM in [51] to formulate a coupled μPPM paradigm for dynamic shear banding and fracturing in multiphase porous media. In this new formulation, it is assumed that the solid skeleton is a micro-polar material while the fluid phases are non-polar following the classical micro-polar poromechanics (e.g., [43, 44]).

---

[*]Corresponding author
*Email address:* xysong@ufl.edu (Xiaoyu Song)

An attractive and salient feature of PPM is that the multiphase discontinuities (i.e., discontinuities in displacement and fluid pressure) can develop naturally based on field equations and material models. For instance, fracture forms naturally when sufficient bonds break at a material point. Here in PPM, the bond means the poromechanics interactions in porous media, which is different from its original definition in PD for solids. An energy-based bond-breakage criterion based on the effective force concept has been formulated to detect bond breakage in PPM (e.g., [29, 32]). Meanwhile, the J-integral can be used in PD for modeling crack propagation in solids [52–54]. In the present study, as the new contributions, we reformulate the J-integral and the energy-based bond breakage criteria incorporating the micro-rotation of material points. As an instability problem, the formation of shear bands in porous media can be detected by the second-order work for multiphase porous media [28]. The classical second-work criterion [55, 56] has been used in PPM through the multiphase correspondence principle to detect shear bands in porous media [29]. In the present study, a nonlocal second-order work for the μPPM paradigm is formulated for the first time directly using the effective force state and moment state [27]. Furthermore, the stabilized Cosserat PPM correspondence principle [51] is used to incorporate the classical micro-polar constitutive models for the solid skeleton, and the stabilized multiphase PPM correspondence principle [30] is used to modeling fluid flow through the classical non-polar fluid flow model.

We have numerically implemented the novel μPPM paradigm through an explicit two-way fractional-step algorithm in time [29] and a hybrid Lagrangian-Eulerian meshfree method in space with Open MPI [57] for high-performance computing. It is noted that the two-way fractional-step algorithm splits the coupled problem into a deformation/fracture problem and an unsaturated fluid flow problem in parallel. We refer to the celebrated literature on fractional-step/staggered and monolithic algorithms for numerically implementing the coupled poromechanics in time (e.g., [58–62], among others). Representative numerical examples are presented to validate the implemented μPPM paradigm and demonstrate its efficacy and robustness in modeling dynamic failure and fracturing in unsaturated porous media. As a brief summary, the original contribution of the present study consists of (i) the formulation of a multiphase micro-polar periporomechanics assuming a micro-polar solid skeleton and a non-polar fluid, (ii) the formulation of a nonlocal micro-polar second-work for detecting shear bands and nonlocal micro-polar energy-based and J-integral for fracture propagation in porous media, and (iii) the computational implementation through an explicit parallel staggered meshfree algorithm and validation of the proposed μPPM paradigm.

The outline of this article is as follows. Section 2 deals with the mathematical formulation of the μPPM paradigm that includes the governing equations, the micro-polar J-integral, the energy-based bond breakage criterion, the micro-polar second-work, and stabilized coupled Micropolar PPM constitutive correspondence principle. Section 3 presents the numerical implementation of the proposed μPPM paradigm through an explicit double-way staggered algorithm. Section 4 presents numerical examples to validate the implemented μPPM paradigm and demonstrate its efficacy and robustness in modeling dynamic shear banding and fracturing in unsaturated porous media, followed by a closure in Section 5. For sign convention, the assumption in continuum mechanics is adopted, i.e., the tensile force and deformation under tension are positive, and for pore fluid pressure, compression is positive, and tension is negative.

## 2. Mathematical formulation

This section presents the mathematical formulation of the proposed μPPM. It consists of four parts. Part I deals with the balance equations for the fully coupled fracturing unsaturated μPPM by assuming passive pore pressure (i.e., zero pore air pressure [1, 2]). Part II presents the fracture and failure criteria, including the micro-polar J-integral, energy-based bond breakage criteria, and the second-order work considering the micro-rotational degree of freedom. Part III deals with the μPPM constitutive correspondence principle with stabilization through which the classical micro-polar material models can be incorporated into the μPPM paradigm. Part IV introduces the classical micro-polar elastic and plastic models for the solid skeleton.

### 2.1. Balance equations for the fully coupled fracturing unsaturated μPPM paradigm

This part concerns the balance equations of the proposed unsaturated micro-periporomechanics paradigm. In line with non-polar periporomechanics, the multiphase porous media is assumed to be represented by a finite number of mixed material points (i.e., solid and fluid material points superimposed). The kinematics of the mixed materials are described in the Relative-Eulerian Lagrangian framework, i.e., the solid material points in Lagrangian and the fluid material points in the Eulerian relative to the Lagrangian for solid material points. The μPPM paradigm postulates that the solid skeleton is micro-polar and the fluid phase is non-polar. Thus, each mixed material point has three types of degrees of freedom, i.e., displacement, micro-rotation, and pore fluid pressures. Following the classical poromechanics, the density of the unsaturated porous media assuming weightless pore air can be written as

$$\rho = (1 - \phi)\rho_s + S_r \phi \rho_w, \tag{1}$$

where $\rho_s$ is the intrinsic density of solid phase, $\rho_w$ in the intrinsic density of fluid phase, $\phi$ is the porosity, and $S_r$ is the degree of saturation. In what follows, we first present the kinematics of an unsaturated material through the μPPM paradigm. Second, we present the balance equations for material points in the bulk (i.e., non-fracturing zone). Third, we present the balance equations for material pints in the fracture zone.

#### 2.1.1. Kinematics of unsaturated μPPM materials

In μPPM, it is postulated that a material point interacts with material points within its family. In this study, the family is assumed as a sphere centered at the material point under consideration. The sphere's radius is called the horizon, δ, a positive number. In other words, it is assumed that two material points at a finite distance have poromechanical interactions, and the horizon is the maximum distance of such two material points

Let x and x' be the position vectors of two solid material points x and x', respectively, in the reference configuration. Let y and y' be the position vectors of the two solid material points in the deformed configuration, respectively. For notation simplicity, in the remaining presentation, a variable with no prime is associated with material point x while a variable with a prime is associated with material point x'. The reference position state of the bond ξ between x and x' in the reference configuration is written as

$$\underline{X} = x' - x. \tag{2}$$

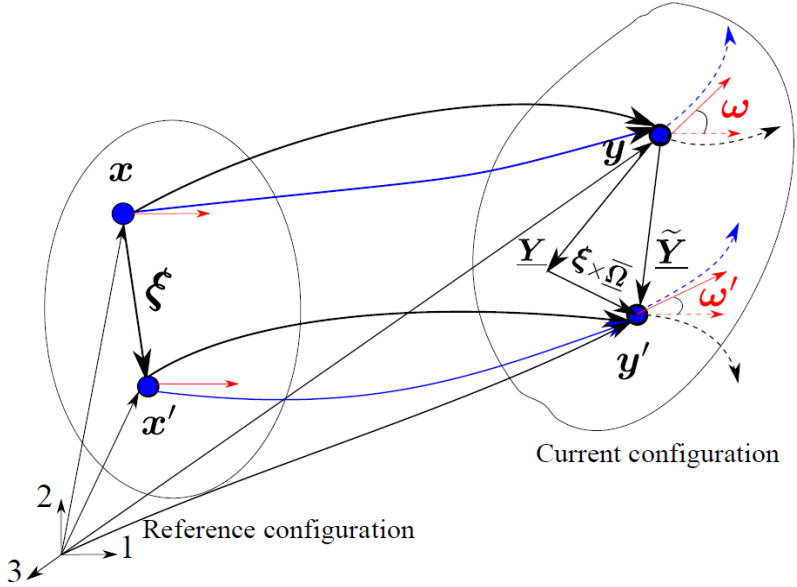

Figure 1: Kinematics of the solid skeleton shown in black and pore fluid material points shown in blue

The deformation state and the relative displacement state in the current configuration read

$$\underline{Y} = y' - y, \tag{3}$$
$$\underline{U} = u' - u. \tag{4}$$

Let $\widehat{\omega}$ and $\widehat{\omega}'$ be the micro-rotations of the two solid material points in the current configuration. Assuming null initial micro-rotations, the relative micro-ration state reads

$$\underline{\Omega} = \widehat{\omega}' - \widehat{\omega} \tag{5}$$

let $\overline{\Omega}$ be the average rotation state that is defined

$$\overline{\Omega} = \frac{1}{2}\left(\widehat{\omega}' + \widehat{\omega}\right). \tag{6}$$

The composite displacement state can be defined as

$$\widehat{\underline{\mathscr{U}}} = \underline{\mathscr{U}} - \overline{\Omega} \times \underline{X} \tag{7}$$

The fluid pressure state can be defined as

$$\underline{\Phi} = p'_w - p_w. \tag{8}$$

Let $p_w$ and $p'_w$ be the water pressure of the two solid material points in the current configuration. Next, we present the field equations for the coupled μPPM.

### 2.1.2. Governing equations in the bulk space

Through the effective force state concept [27] and following [51], the equations of motion and the moment balance for unsaturated micro-polar porous media can be written as

$$\rho \ddot{u} = \int_{\mathcal{H}} \left(\underline{\overline{T}} - \underline{\overline{T}}'\right) d\mathcal{V}' - \int_{\mathcal{H}} \left(S_r \underline{T}_w - S'_r \underline{T}'_w\right) d\mathcal{V}' + \rho g, \quad (9)$$

$$\mathcal{I}^s \ddot{\hat{\omega}} = \int_{\mathcal{H}} \left(\underline{M} - \underline{M}'\right) d\mathcal{V}' + \frac{1}{2} \int_{\mathcal{H}} \underline{Y} \times \left[\left(\underline{\overline{T}} - S_r \underline{T}_w\right) - \left(\underline{\overline{T}}' - S'_r \underline{T}'_w\right)\right] d\mathcal{V}' + l, \quad (10)$$

where $\ddot{u}$ is acceleration, $\underline{\overline{T}}$ is the effective force state, $\underline{T}_w$ is the fluid force state, $g$ is gravitation acceleration, $I_s$ is micro-inertia of solid phase, $\ddot{\hat{\omega}}$ is angular acceleration, $\underline{M}$ is the moment state, $l$ refer to body couple density, and $S_r$ is the degree of saturation, which can be determined from the soil-water retention curve. In this study, the soil-water retention curve of unsaturated porous media assuming passive pore air pressure (e.g., [2, 63–67]) is written as

$$S_r = \left[1 + \left(\frac{-p_w}{s_a}\right)^n\right]^{-m}, \quad (11)$$

where $s_a$ is the scaling factor, m and n are fitting parameters. We refer to [51] for deriving equations (9) and (10) from a free energy function for dry micro-polar porous media.

In geomechanics, $I_s \ddot{\hat{\omega}}$ is the angular momentum of a spinning grain. By considering the mean grain size in porous media as the Cosserat length scale $l$ [41], the micro inertia can be written as

$$\mathcal{I}^s = \frac{\pi}{2}(1 - \phi)\rho_s l^2. \quad (12)$$

Assuming an incompressible solid phase and a non-polar fluid phase, the mass balance equation of the mixture [27] is written as

$$\phi \frac{dS_r}{dt} + S_r \dot{\mathcal{V}}_s + \frac{1}{\rho_w} \int_{\mathcal{H}} \left(\underline{\mathcal{Q}} - \underline{\mathcal{Q}}'\right) d\mathcal{V}' + \mathcal{Q}_s = 0, \quad (13)$$

where $\dot{V}_s$ is the solid volume change rate, Q and Q' are the fluid flow states at x and x', respectively, and $Q_s$ is a source term. It should be noted that the micro-rotation of the solid phase does not affect the volume change rate of the solid phase [44].

Next, we cast the governing equations for material points in the fracturing zone.

### 2.1.3. Governing equations for material points in the fracturing zone

Following (9) and (10), we can write out the equation of motion and the moment balance equation for a fracturing material point as

$$\rho \ddot{u} = \int_{\mathcal{H}} \left(\underline{\overline{T}} - \underline{\overline{T}}'\right) d\mathcal{V}' - \int_{\mathcal{H}} \left(S_l \underline{T}_l - S'_l \underline{T}'_l\right) d\mathcal{V}' + \rho g, \quad (14)$$

$$\mathcal{I}^s \ddot{\hat{\omega}} = \int_{\mathcal{H}} \left(\underline{M} - \underline{M}'\right) d\mathcal{V}' + \frac{1}{2} \int_{\mathcal{H}} \underline{Y} \times \left[\left(\underline{\overline{T}} - S_l \underline{T}_l\right) - \left(\underline{\overline{T}}' - S'_l \underline{T}'_l\right)\right] d\mathcal{V}' + l, \quad (15)$$

where

$$S_l \underline{T}_l = \begin{cases} S_{r,f} \underline{T}_f, & \text{if } D > D_{cr} \ \& \ D' > D_{cr}, \\ S_r \underline{T}_w, & \text{otherwise}, \end{cases} \quad (16)$$

and $S_{r,f}$ is the degree of saturation of the fracturing mixed material point, which can be determined from the soil-water retention curve, D is the damage state variable, and $D_{cr}$ is the critical value of the damage state variable.

Assuming a non-polar fluid phase, the mass balance equation in the fracturing space [27] can be written as

$$\frac{\partial S_{r,f}}{\partial t} + \frac{1}{\rho_w} \int_{\mathcal{H}} \left(\underline{\mathcal{Q}}_f - \underline{\mathcal{Q}}'_f\right) d\mathcal{V}' - \mathcal{Q}_s = 0, \tag{17}$$

where $Q_f$ is the fluid flow state at a fracturing mixed material point.

The fluid flow from the fracturing pore space to the continuous pore space can be determined, assuming that the fluid flows from the pore space into the fracture space. It is assumed that the direction of unsaturated fluid flow is normal to the fracture surface. Thus, following the generalized Darcy's law [29], $Q_s$ can be written as

$$\mathcal{Q}_s = -\frac{k^r k_w}{\mu_w} \left(\frac{p_f - p_w}{l_x}\right) \mathcal{A}/\mathcal{V}, \tag{18}$$

where $p_w$ is pore water pressure in the bulk and pf is pore water pressure in the fracture space, $A$ is the cross-sectional area and $V$ is the volume of a material point, $lx = d/2$, and $d$ is the dimension of a cubic material point [29]. In (18), it is assumed that porosity $\phi = 0$ and the volume coupling term vanishes in the fracture pore space.

## 2.2. μPPM fracturing and instability criteria for unsaturated porous media

In this subsection, in the μPPM framework, we reformulate the micro-polar J integral and energy-based bond breakage criterion for fracturing and the micro-polar second-order work criterion for instability in porous media.

### 2.2.1. μPPM J integral for unsaturated porous media

In this part, we present the micro-polar J integral in the coupled μPPM paradigm following the lines in [35, 45, 47]. The effective force state concept is used to express the force state on the solid skeleton. Figure 2 shows a subregion $B_1$ in the reference configuration $B$. As shown in Figure 2, for a bounded subregion $B_1$ with a constant shape moves with velocity $\dot{u}$ in the x direction (parallel to crack line) in the reference configuration $B$, there is a flux of material through the boundary $\partial B_1$. By assuming the steady-state motion [52], the deformation position of each point and micro rotation of each point inside the subregion $B_1$ is written as

$$\boldsymbol{y} = \boldsymbol{x} + \boldsymbol{u}\Big|_{\boldsymbol{x} - \dot{\bar{\boldsymbol{u}}}t}, \tag{19}$$

$$\widehat{\omega} = \widehat{\omega}\Big|_{\boldsymbol{x} - \dot{\bar{\boldsymbol{u}}}t}. \tag{20}$$

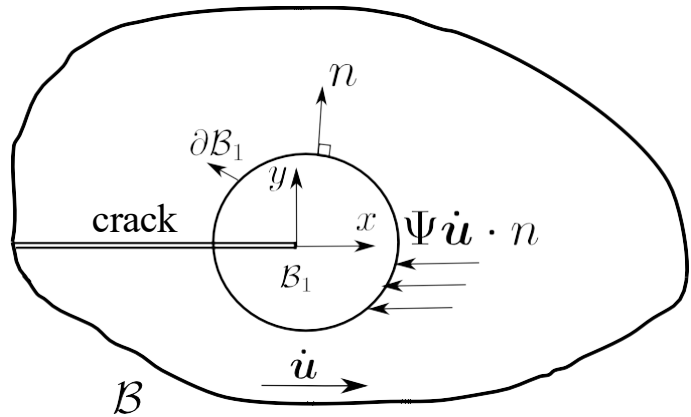

Figure 2: Subregion $B_1$ with constant shape around the crack tip moves with velocity $\boldsymbol{u}^{\cdot}$ in $x$ direction and a flux of material $\psi \boldsymbol{u}^{\cdot} \cdot \boldsymbol{n}$ through the boundary $\partial B_1$

For simplicity, we neglect the body force, kinetic energy, heat transport, and heat sources and assume an isothermal process. Thus, the rate form of the free energy for the solid phase is written as

$$\dot{\Psi} = \Psi_{\underline{\widetilde{\mathcal{Y}}}} \bullet \dot{\underline{\widetilde{\mathcal{Y}}}} + \Psi_{\underline{\Omega}} \bullet \underline{\dot{\Omega}} + \Psi_{\underline{D}} \bullet \underline{\dot{D}}, \tag{21}$$

where $\bullet$ is the dot product of two vector states [35], $\underline{D}$ is the damage scalar state as mentioned previously that cannot decrease over time, and $\Psi_{(\cdot)}$ is the Fŕechet derivative of $\Psi$ to its state variables [35]. It follows from (21) and assuming an elastic deformation that the effective force state and the moment state [27, 51] can be defined as

$$\underline{\overline{\mathcal{T}}} = \Psi_{\underline{\widetilde{\mathcal{Y}}}}, \tag{22}$$

$$\underline{\mathcal{M}} = \Psi_{\underline{\Omega}}. \tag{23}$$

Substituting (22) and (23) into (21), the rate form of the free energy reads

$$\dot{\Psi} = \underline{\overline{\mathcal{T}}} \bullet \left( \underline{\dot{\mathcal{Y}}} - \underline{\dot{\Omega}} \times \underline{\xi} \right) + \underline{\mathcal{M}} \bullet \underline{\dot{\Omega}} + \Psi_{\underline{D}} \bullet \underline{\dot{D}}. \tag{24}$$

We note that by considering the steady state motion, the total time derivative of each quantity over B1 vanishes. Thus, applying the nonlocal Reynolds transport theorem [52] to the free energy
we have

$$\frac{d}{dt} \int_{\mathcal{B}_1} \Psi d\mathcal{V} = \int_{\mathcal{B}_1} \dot{\Psi} d\mathcal{V} + \int_{\partial \mathcal{B}_1} \Psi \dot{u} \cdot n dA = 0. \tag{25}$$

where $\dot{u}$ is the constant velocity, n is the outward unit normal to the boundary $\partial B1$, and A is the area of the boundary $\partial B1$ (see Figure 2). With (24), the first term on the right-hand side of (25) can be written as

$$\int_{\mathcal{B}_1} \dot{\Psi} d\mathcal{V} = \int_{\mathcal{B}_1} \int_{\mathcal{B}/\mathcal{B}_1} \left( \underline{\overline{\mathcal{T}}} \cdot \dot{y}' - \underline{\overline{\mathcal{T}}}' \cdot \dot{y} \right) d\mathcal{V}' d\mathcal{V}$$
$$+ \frac{1}{2} \int_{\mathcal{B}_1} \int_{\mathcal{B}/\mathcal{B}_1} \left[ \underline{\xi} \times \left( \underline{\overline{\mathcal{T}}} \cdot \dot{\widehat{\omega}}' - \underline{\overline{\mathcal{T}}}' \cdot \dot{\widehat{\omega}} \right) \right] d\mathcal{V}' d\mathcal{V}$$
$$+ \int_{\mathcal{B}_1} \int_{\mathcal{B}/\mathcal{B}_1} \left( \underline{\mathcal{M}} \cdot \dot{\widehat{\omega}}' - \underline{\mathcal{M}}' \cdot \dot{\widehat{\omega}} \right) d\mathcal{V}' d\mathcal{V} + \int_{\mathcal{B}_1} \left( \Psi_{\underline{D}} \bullet \underline{\dot{D}} \right) d\mathcal{V}. \tag{26}$$

Substituting (26) into (25), we obtain the expression for the damage energy dissipation as

$$\int_{\mathcal{B}_1} \left( \Psi_{\underline{D}} \bullet \underline{\dot{D}} \right) d\mathcal{V} = \int_{\mathcal{B}_1} \int_{\mathcal{B}/\mathcal{B}_1} \left[ \left( \underline{\overline{\mathcal{T}}} \cdot u'_x - \underline{\overline{\mathcal{T}}}' \cdot u_x \right) \cdot \dot{u} \right] d\mathcal{V}' d\mathcal{V}$$
$$+ \frac{1}{2} \int_{\mathcal{B}_1} \int_{\mathcal{B}/\mathcal{B}_1} \left[ \underline{\xi} \times \left( \underline{\overline{\mathcal{T}}} \cdot \widehat{\omega}'_x - \underline{\overline{\mathcal{T}}}' \cdot \widehat{\omega}_x \right) \cdot \dot{u} \right] d\mathcal{V}' d\mathcal{V}$$
$$+ \int_{\mathcal{B}_1} \int_{\mathcal{B}/\mathcal{B}_1} \left[ \left( \underline{\mathcal{M}} \cdot \widehat{\omega}'_x - \underline{\mathcal{M}}' \cdot \widehat{\omega}_x \right) \cdot \dot{u} \right] d\mathcal{V}' d\mathcal{V} - \int_{\partial \mathcal{B}_1} \left( \Psi \overline{u} \cdot n \right) dA. \tag{27}$$

where $(.)x$ is the directional derivative in the direction of the stead state motion. Here, it is assumed that the steady flow is in the x direction for simplicity. The relation between the integral and the rate of energy dissipation [52] can be written as

$$\mathcal{J} \cdot \dot{\overline{u}} = - \int_{\mathcal{B}_1} \Psi_{\underline{D}} \bullet \underline{\dot{D}} d\mathcal{V}. \tag{28}$$

It follows from (27) and (28) that the **J** integral for the μPPM paradigm reads

$$\mathcal{J} = \int_{\partial \mathcal{B}_1} (\Psi \boldsymbol{n}) \, d\mathcal{A} - \int_{\mathcal{B}_1} \int_{\mathcal{B}/\mathcal{B}_1} \left( \underline{\mathcal{T}} \cdot \boldsymbol{u}'_x - \underline{\mathcal{T}}' \cdot \boldsymbol{u}_x \right) d\mathcal{V}' d\mathcal{V}$$

$$- \frac{1}{2} \int_{\mathcal{B}_1} \int_{\mathcal{B}/\mathcal{B}_1} \left[ \boldsymbol{\xi} \times \left( \underline{\mathcal{T}} \cdot \widehat{\boldsymbol{\omega}}'_x - \underline{\mathcal{T}}' \cdot \widehat{\boldsymbol{\omega}}_x \right) \right] d\mathcal{V}' d\mathcal{V}$$

$$- \int_{\mathcal{B}_1} \int_{\mathcal{B}/\mathcal{B}_1} \left( \underline{\mathcal{M}} \cdot \widehat{\boldsymbol{\omega}}'_x - \underline{\mathcal{M}}' \cdot \widehat{\boldsymbol{\omega}}_x \right) d\mathcal{V}' d\mathcal{V}. \tag{29}$$

We note that (29) is the μPPM equivalent of the J integral in the standard micro-polar continuum theory. different from the J integral for peridynamic solids [52], the μPPM J integral incorporates the micro-rotations of material points. We note that u and $\widehat{\omega}$ are differentiable only in the directions of motion and rotation there are used. Thus, it is permissible to have discontinuities parallel to the direction of crack propagation. Furthermore, J is path independent" in the sense that B1 under deformation can include any number of additional material points in which there is no energy dissipation (i.e., the value of J not changing). The energy dissipation per unit crack area can be determined from [52]. For simplicity, let the body be a plate with a crack through its thickness that propagates along the x direction. The energy dissipated per unit crack area can be written as

$$\mathcal{G} = \mathcal{J} \cdot \hat{n}/h, \tag{30}$$

where $\hat{n}$ is the unit vector parallel to the crack surface and h is the plate thickness. Substituting (29) into (30) generate

$$\mathcal{G} = \int_{\partial \mathcal{B}_1} (\Psi \boldsymbol{n} \cdot \hat{n}) \, d\mathcal{A}/h - \int_{\mathcal{B}_1} \int_{\mathcal{B}/\mathcal{B}_1} \left[ \left( \underline{\mathcal{T}} \cdot \boldsymbol{u}'_x - \underline{\mathcal{T}}' \cdot \boldsymbol{u}_x \right) \cdot \hat{n} \right] d\mathcal{V}' d\mathcal{V}/h$$

$$- \frac{1}{2} \int_{\mathcal{B}_1} \int_{\mathcal{B}/\mathcal{B}_1} \left[ \boldsymbol{\xi} \times \left( \underline{\mathcal{T}} \cdot \widehat{\boldsymbol{\omega}}'_x - \underline{\mathcal{T}}' \cdot \widehat{\boldsymbol{\omega}}_x \right) \cdot \hat{n} \right] d\mathcal{V}' d\mathcal{V}/h \tag{31}$$

$$- \int_{\mathcal{B}_1} \int_{\mathcal{B}/\mathcal{B}_1} \left[ \left( \underline{\mathcal{M}} \cdot \widehat{\boldsymbol{\omega}}'_x - \underline{\mathcal{M}}' \cdot \widehat{\boldsymbol{\omega}}_x \right) \cdot \hat{n} \right] d\mathcal{V}' d\mathcal{V}/h. \tag{32}$$

It is noted that from the linear elastic fracture mechanics, the critical energy release rate per unit crack area Gcr for mode I crack reads

$$\mathcal{G}_{cr} = \mathcal{K}_I^2 \frac{1 - \nu^2}{E}, \tag{33}$$

where $k_I$ is the fracture toughness of mode I crack, E is Young's modulus and $\nu$ is Poisson ratio. Next, we present the energy-based μPPM bond breakage criterion for micro-polar unsaturated porous media.

*2.2.2. Energy-based μPPM bond breakage criterion*

This part presents an energy-based bond breakage criterion breakage through the effective force state concept for unsaturated porous media [27]. The bond-breakage criterion depends on the deformation energy in a bond. The effective force and moment states that are the energy conjugates to the composite state and relative rotation state, respectively, are used to determine the deformation energy of a bond [29]. Thus, the energy density in the bond of porous media can be obtained as

$$\mathcal{W} = \int_0^t \left[ \left( \underline{\mathcal{T}} + S_r \underline{\mathcal{T}}_w \right) - \left( \underline{\mathcal{T}}' + S_r' \underline{\mathcal{T}}'_w \right) \right] \dot{\widehat{\mathcal{U}}} \, dt$$

$$+ \int_0^t \left( \underline{\mathcal{M}} - \underline{\mathcal{M}}' \right) \dot{\boldsymbol{\Omega}} \, dt. \tag{34}$$

where *t* is the loading time. The bond breakage is modeled through the influence function at the material point level. The influence function of the broken bond will be updated by $\underline{\varrho}\omega$, where $\underline{\varrho}$ is determined by

$$\underline{\varrho} = \begin{cases} 0 & \text{if } \mathcal{W} \geqslant \mathcal{W}_{cr}, \\ 1 & \text{otherwise} \end{cases} \qquad (35)$$

where $\mathcal{W}_{cr}$ is the critical bond energy density. The critical energy density for bond breakage can be determined from the critical energy release rate as

$$\mathcal{W}_{cr} = \frac{4\mathcal{G}_{cr}}{\pi \delta^4}. \qquad (36)$$

In μPPM, the damage of the solid skeleton material point is tracked through a local scalar damage variable $\varphi$. This damage variable is defined as the fraction of broken bonds at a material point

$$\varphi = 1 - \frac{\int_{\mathcal{H}} \underline{\varrho}\omega d\mathcal{V}'}{\int_{\mathcal{H}} \underline{\omega} d\mathcal{V}'}. \qquad (37)$$

As in [29], the space between material points x and x' is assumed fractured where $w > w_{cr}$ and D>Dcr and D'>Dcr hold for both material points. The two material points are defined as fracture mixed points. The fractured mixed material points have bulk fluid pressure and fracture fluid pressure. These fluid pressure are utilized to model unsaturated fluid flow in fractured space. Figure 3 plots the kinematics of a bond across fracture space. As shown in Figure 3, the micro- polar relative displacement vector $\eta\overline{\Omega}$ is decomposed into two components [29, 68] where the first component represents the opening displacement and the second one is the dislocation of the crack. The crack aperture c, which is related to the opening displacement, can be written as

$$\underline{c} = \underline{\mathcal{Y}}_{\overline{\Omega}} \cos\beta - \underline{x}, \qquad (38)$$

where angle $\beta$ is shown in Figure 3. Therefore, the crack width at fracture point x can be approximated by the average of bond apertures of all broken bonds as

$$a_f = \frac{\int_{\mathcal{H}} \omega \underline{c} d\mathcal{V}'}{\int_{\mathcal{H}} \underline{\omega} d\mathcal{V}'}, \qquad (39)$$

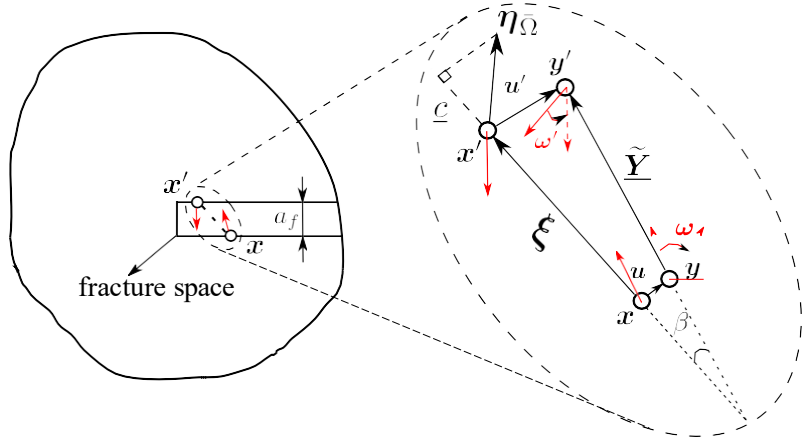

Figure 3: Kinematics of a bond across fracture space

*2.2.3. Nonlocal μPPM second-order work*

The second-order work criterion can be used to detect material instability, including shear banding in porous media [55, 56]. It states that the material loses stability if the second-order work is negative. We have formulated the nonlocal second-order work for the micro-polar PPM paradigm through the Cosserat PPM correspondence principle in [51]. The second-order work d2W can be written as

$$d^2\mathcal{W} = d\overline{\boldsymbol{\sigma}} : d\boldsymbol{\varepsilon} + d\boldsymbol{m} : d\boldsymbol{\kappa}. \qquad (40)$$

In this study, we express the μPPM second-order work from the internal energy rate of the skeleton through the effective force and relative moment states. The rate of the internal energy density of the solid skeleton $\mathcal{E}$ [51] can be written as

$$\dot{\mathcal{E}} = \underline{\overline{\mathcal{T}}} \bullet \dot{\widehat{\underline{\mathcal{U}}}} + \underline{\mathcal{M}} \bullet \dot{\underline{\Omega}}. \tag{41}$$

Thus, following the lines in [31] the μPPM second-order work $d\mathcal{E}$ in terms of the effective force state and relative moment state can be written as

$$d^2\mathcal{E} = \int_{\mathcal{H}} \left( \Delta\underline{\overline{\mathcal{T}}} \cdot \Delta\widehat{\underline{\mathcal{U}}} + \Delta\underline{\mathcal{M}} \cdot \Delta\underline{\Omega} \right) d\mathcal{V}'. \tag{42}$$

We note that (42) and (40) are equivalent for the μPPM correspondence material model. Next, we present the stabilized multiphase μPPM correspondence principle.

*2.3. Stabilized multiphase μPPM correspondence principle*

The constitutive models for both the solid phase and fluid phase are needed to close the governing equations presented in Sections 2.1.2 and 2.1.3, e.g., the relationships between the effective force state and the composite deformation state, the relative moment state and the micro-rotation state, and the fluid flow state and the relative fluid pressure state. We note that the new constitutive models for μPPM can be cast following the lines in [27, 35]. However, it would require extensive experimental tests to validate the constitutive models for their practical applications, which is beyond the scope of the present study. In this study, we adopt the micro- polar constitutive correspondence principle in [51] and the multiphase constitutive principle in [30] to incorporate the classical micro-polar constitutive model for the solid skeleton and the classical non-polar fluid flow model for the fluid phase. Next, we introduce the multiphase μPPM correspondence principle, followed by a multiphase stabilization scheme.

*2.3.1. Multiphase μPPM correspondence principle*

This part presents the multiphase μPPM constitutive correspondence principle. First, we determine the effective force and moment states through the Cosserat PPM constitutive principle for the solid phase. We note that the Cosserat PPM constitutive correspondence principle was proposed in [51] by equating the deformation energy of the Cosserat PPM and the classical strain energy of the classical micropolar continuum for the solid skeleton. In this study, for conciseness, we omit the derivation and refer to [51] for the detailed derivation for the solid skeleton. Following the effective force state concept for unsaturated porous media [27], the total force state and the moment state can be written as

$$\underline{\mathcal{T}} = \omega\boldsymbol{\sigma}\mathcal{K}^{-1}\underline{\boldsymbol{\xi}}, \tag{43}$$

$$\underline{\mathcal{M}} = \omega m\mathcal{K}^{-1}\underline{\boldsymbol{\xi}}, \tag{44}$$

where $\sigma$ is the total stress of unsaturated porous media, is the unit weighting function, m is the couple stress tensor, and K is the shape function, which is defined as

$$\mathcal{K} = \int_{\mathcal{H}} \underline{\boldsymbol{\xi}} \otimes \underline{\boldsymbol{\xi}} d\mathcal{V}'. \tag{45}$$

The total stress state of the classical unsaturated poromechanics [2] assuming passive pore air pressure can be written as

$$\boldsymbol{\sigma} = \overline{\boldsymbol{\sigma}} - S_r p_w \mathbf{1}, \tag{46}$$

where $\overline{\sigma}$ is the effective stress and 1 is the second-order identity tensor. Note that in (46) the pore water pressure is negative in unsaturated porous media. In PPM, the total force state of unsaturated porous media [28] can be written as

$$\underline{\mathcal{T}} = \overline{\underline{\mathcal{T}}} + S_r \underline{\mathcal{T}}_w. \tag{47}$$

Then, it follows from (43), (46), and (47) that the effective force state and the fluid force state can be written as

$$\overline{\underline{\mathcal{T}}} = \omega \overline{\boldsymbol{\sigma}} \mathcal{K}^{-1} \underline{\boldsymbol{\xi}}, \tag{48}$$

$$\underline{\mathcal{T}}_w = -\underline{\omega}(p_w \mathbf{1}) \mathcal{K}^{-1} \underline{\boldsymbol{\xi}}. \tag{49}$$

From the classical micro-polar poromechanics, the effective stress and the couple stress tensor can be determined the classical micro-polar constitutive models given the strain tensor " and the wryness tensor (e.g., [42–44, 69]). The nonlocal versions of $\varepsilon$ and $\kappa$ are written as

Next, we determine the fluid flow states in the bulk and fracturing porous space through the constitutive correspondence principle for fluid flow in unsaturated porous media. Following the constitutive correspondence principle for fluid flow [28], the fluid flow states Q and Qf in the bulk space and the fracturing space can be written as

$$\underline{\mathcal{Q}} = \underline{\omega} \left( \rho_w \boldsymbol{q}_w \right) \mathcal{K}^{-1} \underline{\boldsymbol{\xi}}, \tag{50}$$

$$\underline{\mathcal{Q}}_f = \underline{\omega} \left( \rho_w \boldsymbol{q}_f \right) \mathcal{K}^{-1} \underline{\boldsymbol{\xi}}, \tag{51}$$

where $q_w$ and $q_f$ are the fluid flux vectors in the bulk and fracture spaces, respectively. The fluid flow vectors can be obtained from the generalized Darcy's law for unsaturated fluid flow as

$$\boldsymbol{q}_w = -\frac{k_w^r k_w}{\mu_w} \widetilde{\nabla \Phi}, \tag{52}$$

$$\boldsymbol{q}_f = -\frac{k_f^r k_f}{\mu_w} \widetilde{\nabla \Phi}_f, \tag{53}$$

where $k_w$ is the intrinsic permeability of the bulk space, $k_r$ is the relative permeability of the bulk space, µw is the water viscosity, $\widetilde{\nabla \Phi}$ is the non-local fluid pressure gradient of the bulk space, $k_r^f$ is the relative permeability of the fracture space, and kf is the intrinsic permeability of the fracture space, and $\widetilde{\nabla \Phi}_f$ is the fluid pressure gradient in the fracture space. The relative permeabilities are determined by

$$k_w^r = \sqrt{S_r} \left[ 1 - \left( 1 - S_r^{1/m} \right)^m \right]^2, \tag{54}$$

$$k_f^r = S_{r,f}^{1/2} \left[ 1 - \left( 1 - S_{r,f}^{1/m} \right)^m \right]^2, \tag{55}$$

where $S_{r,f}$ is the degree of saturation of fluid phase in fracturing space and m is the same material parameter as defined in (11). The intrinsic permeability of the fracture can be determined by the so-called cubic law [7, 28] as

$$k_f = \frac{a_f^2}{12}, \tag{56}$$

where $a_f$ is the crack aperture.

The nonlocal fluid pressure gradients are written as

$$\widetilde{\nabla \Phi} = \left( \int_{\mathcal{H}} \omega \underline{\Phi} \underline{\boldsymbol{\xi}} d\mathcal{V}' \right) \mathcal{K}^{-1}, \tag{57}$$

$$\widetilde{\nabla \Phi}_f = \left( \int_{\mathcal{H}} \omega \underline{\Phi}_f \underline{\boldsymbol{\xi}} d\mathcal{V}' \right) \mathcal{K}^{-1}, \tag{58}$$

where the fluid pressure state in the fracture space is

$$\underline{\Phi}_f = p'_f - p_f. \tag{59}$$

*2.3.2. Multi-phase stabilization scheme*

The energy method is used to eliminate the zero-energy mode instability in the correspondence constitutive principle [30]. We first present the residuals (i.e., non-uniform parts) of the composite deformation and relative rotation states as

$$\underline{\mathscr{R}}_1 = \underline{\mathscr{U}}_{\overline{\Omega}} - \varepsilon \underline{\xi}, \tag{60}$$

$$\underline{\mathscr{R}}_2 = \underline{\Omega} - \kappa \underline{\xi}, \tag{61}$$

Then, the stabilized force and moment states [30] can be written as

$$\underline{\overline{\mathscr{T}}} = \underline{\omega \xi} \overline{\sigma} \mathscr{K}^{-1} \underline{\xi} + \alpha_1 \underline{\mathscr{R}}_1, \tag{62}$$

$$\underline{\mathscr{M}} = \underline{\omega \xi} m \mathscr{K}^{-1} \underline{\xi} + \alpha_2 \underline{\mathscr{R}}_2. \tag{63}$$

where $\alpha_1$ and $\alpha_2$ are two scalar states. The two scalar states can be written as

$$\alpha_1 = \frac{g\mathscr{C}_1}{\omega_0} \underline{\omega}, \tag{64}$$

$$\alpha_2 = \frac{g\mathscr{C}_2}{\omega_0} \underline{\omega}, \tag{65}$$

where

$$\omega_0 = \int_{\mathscr{H}} \underline{\omega \xi} d\mathscr{V}', \tag{66}$$

$g$ is a number between 0 and 1. C1 and C2 are two material parameters. For three dimensional case, these two material parameters are written as

$$\mathscr{C}_1 = \frac{12\mathscr{D}}{|\underline{\xi}|^3}, \tag{67}$$

$$\mathscr{C}_2 = \frac{\mathscr{D}}{|\underline{\xi}|}. \tag{68}$$

where D is a parameter that depends on the length scale. In three-dimensions it reads

$$\mathscr{D} = \frac{E(1-4\nu)}{4\pi \delta^2 (1-\nu-2\nu^2)}. \tag{69}$$

Similarly, for the fluid phase, we define the residuals (i.e., the non-uniform parts) of the fluid pressure states as

$$\underline{\mathscr{R}}_w = \underline{\Phi} - \widetilde{\nabla \underline{\Phi}} \underline{\xi}, \tag{70}$$

$$\underline{\mathscr{R}}_f = \underline{\Phi}_f - \widetilde{\nabla \underline{\Phi}_f} \underline{\xi}, \tag{71}$$

Then, the fluid flow state with stabilized can be written as

$$\underline{\mathscr{Q}} = \underline{\omega} \left( q_w \mathscr{K}^{-1} \underline{\xi} + \alpha_3 \underline{\mathscr{R}}_w \right), \tag{72}$$

$$\underline{\mathscr{Q}}_f = \underline{\omega} \left( q_f \mathscr{K}^{-1} \underline{\xi} + \alpha_4 \underline{\mathscr{R}}_f \right), \tag{73}$$

In (72) and (73),

$$\alpha_3 = \frac{g\mathscr{C}_3}{\omega_0} \underline{\omega}, \tag{74}$$

$$\alpha_4 = \frac{g\mathscr{C}_4}{\omega_0} \underline{\omega}, \tag{75}$$

where $C_3$ and $C_4$ are the micro-conductivities [30]. For the three-dimensional case, they are written as

$$\mathcal{C}_3 = \frac{6k_w}{\pi \delta^4}, \tag{76}$$

$$\mathcal{C}_4 = \frac{6k_f}{\pi \delta^4}. \tag{77}$$

Next, we present the classical micro-polar elastic and plastic constitutive models adopted in this study.

### 2.4. Classical micro-polar material models
#### 2.4.1. Micropolar linear elastic model

For the micropolar elastic model for the porous media, the effective stress tensor $\bar{\sigma}_{ij}$ and the couple stress tensor $m_{ij}$ [39] can be written as

$$\bar{\sigma}_{ij} = \lambda \varepsilon^e_{kk} + (\mu + \mu_c)\varepsilon^e_{ij} + (\mu - \mu_c)\varepsilon^e_{ji}, \tag{78}$$

$$m_{ij} = \frac{1}{2}\mu l^2 \kappa^e_{ij}, \tag{79}$$

where i, j, k = 1, 2, 3, $\varepsilon_{ij}$ is the elastic strain tensor, $\kappa_{ij}$ is the elastic wryness tensors, $\lambda$ is L'ame's first elastic constant, $\mu$ is the shear modulus, $\mu c$ is the micro-polar shear modulus and $l$ is the micro-polar length scale. We note that the micro-polar length scale can depend on the porous media's micro-structure. How to determine the micro-polar length scale from the microstructure of granular materials is beyond the scope of the present study. We refer to the distinguished literature on the subject (e.g., [40, 41]).

#### 2.4.2. Micropolar plastic model

In this section, the micropolar plastic model with a Drucker-Prager type yield criterion [42, 69]. The strain tensor $\varepsilon_{ij}$ and the wryness tensor $\kappa_{ij}$ are additively decomposed into elastic and plastic parts as

$$\varepsilon_{ij} = \varepsilon^e_{ij} + \varepsilon^p_{ij}, \tag{80}$$

$$\kappa_{ij} = \kappa^e_{ij} + \kappa^p_{ij}, \tag{81}$$

where $\varepsilon^e_{ij}$ is the elastic strain tensor, $\varepsilon^p_{ij}$ is the plastic strain tensor, $\kappa^e_{ij}$ is the elastic wryness tensor, and $\kappa^e_{ij}$ is the plastic wryness tensors. The effective stress tensor $\bar{\sigma}_{ij}$ and couple stress tensor $m_{ij}$ can be calculated from elastic components of strain and wryness tensors.

$$q = \sqrt{\frac{1}{2}\left(\bar{\sigma}_{ij}\tilde{\sigma}_{ij} + \frac{3m_{ij}m_{ij}}{l^2}\right)}, \tag{82}$$

$$\bar{p} = \frac{1}{3}\bar{\sigma}_{kk}, \tag{83}$$

where

$$\tilde{\sigma}_{ij} = -\bar{\sigma}_{kk} + \frac{3}{2}\left(\bar{\sigma}_{ij} + \bar{\sigma}_{ji}\right). \tag{84}$$

The yield function and the plastic potential are written as

$$f = q + \sqrt{3}a_1 \bar{p} + a_2, \tag{85}$$

$$g = q + \sqrt{3}a_3 \bar{p} + a_2, \tag{86}$$

where

$$a_1 = \frac{2\sin\varphi}{\sqrt{3}(3-\sin\varphi)}, \tag{87}$$

$$a_2 = -\frac{6c\cos\varphi}{\sqrt{3}(3-\sin\varphi)}, \tag{88}$$

$$a_3 = \frac{2\sin\psi}{\sqrt{3}(3-\sin\psi)}, \tag{89}$$

and $\varphi$ is the frictional angle, $\psi$ is the dilatancy angle, and c is cohesion. The parameter c is defined as

$$c = c_0 + h\bar{\varepsilon}^p, \tag{90}$$

where h is the linear isotropic hardening modulus and $\bar{\varepsilon}^p$ is a plastic variable. The non-associative plastic flow rule can be written as

$$\dot{\varepsilon}^p_{ij} = \dot{\lambda}\frac{\partial g}{\partial \bar{\sigma}_{ij}}, \tag{91}$$

$$\dot{\kappa}^p_{ij} = \dot{\lambda}\frac{\partial g}{\partial m_{ij}}. \tag{92}$$

where $\dot{\lambda}$ is the plastic multiplier, which is obtained from consistency condition $\dot{f} = 0$. The internal plastic variable [70] is defined as

$$\dot{\bar{\varepsilon}}^p = \sqrt{\frac{1}{3}\dot{\varepsilon}^p_s : \dot{\varepsilon}^p_s + \frac{1}{3}\dot{\varepsilon}^p_s : \dot{\varepsilon}^{p,T}_s + \frac{2}{3}\dot{\kappa}^p : \dot{\kappa}^p}, \tag{93}$$

where $\varepsilon^p_s$ is the deviatoric part of the plastic strain tensor. In the following section, we present the numerical implementation of the coupled μPPM paradigm.

## 3. Numerical implementation

This section deals with the numerical implementation of the coupled μPPM paradigm through a hybrid Lagrangian-Eulerian meshfree method in space and an explicit-explicit dual-direction fractional-step algorithm in time [29]. The dual-direction fractional-step algorithm splits the fully coupled poromechanics problem into a deformation/fracture problem (i.e., the solid solver) and an unsaturated fluid flow problem (i.e., the fluid solver) in parallel. We refer to the distinguished literature (e.g., [58–62], among others) on fractional-step/staggered and monolithic algorithms for numerically implementing coupled poromechanics problems in time. Here, the dual-direction staggered algorithm means that the solid solver or the fluid solver can be called first, given the coupled poromechanics problem. For instance, for a fluid-driven cracking problem, the program will call the unsaturated fluid flow solver first, assuming no deformation, and then the fluid pressure will be passed to the deformation solver.

### 3.1. Spatial discretization

The coupled governing equations are spatially discretized by the hybrid Lagrangian-Eulerian meshfree scheme, i.e., Lagrangian for the solid phase and Eulerian (relative to the solid) for the fluid phase. Consistent with the mathematical formulations in Section 2, a porous material body is discretized into a finite number of mixed material points (i.e., the solid skeleton and the pore fluid). Each mixed material point is endowed with three types of degrees of freedom, i.e., displacement, rotation, and fluid pressure. The boundary conditions are imposed through the boundary layer method. Let P be the number of total mixed material points in the problem domain and Ni be the number of material points in the horizon of material point i. The spatial discretization of the motion equation and the balance of moment are written as

$$0 = \mathcal{A}_{i=1}^{\mathcal{P}} \left( \mathcal{M}_i \ddot{\mathbf{u}}_i - \overline{\mathcal{T}}_i + \mathcal{T}_{w,i} - \mathcal{M}_i \mathbf{g} \right), \tag{94}$$

$$0 = \mathcal{A}_{i=1}^{\mathcal{P}} \left( \mathcal{L}_i \ddot{\boldsymbol{\omega}}_i - \mathcal{M}_i - \overline{\mathcal{M}}_i - l_i \mathcal{V}_i \right), \tag{95}$$

where

$$\mathcal{M}_i = [\rho_s (1 - \phi_i) + \rho_w S_{r,i} \phi_i] \mathcal{V}_i. \tag{96}$$

$$\overline{\mathcal{T}}_i = \sum_{j=1}^{\mathcal{N}_i} \left( \underline{\mathcal{T}}_{ij} - \underline{\mathcal{T}}_{ji} \right) \mathcal{V}_j \mathcal{V}_i, \tag{97}$$

$$\mathcal{T}_{w,i} = \sum_{j=1}^{\mathcal{N}_i} (1 - \Gamma_{ij}) \left( S_{r,i} \underline{\mathcal{T}}_{w,ij} - S_{r,j} \underline{\mathcal{T}}_{w,ji} \right) \mathcal{V}_j \mathcal{V}_i$$

$$+ \sum_{j=1}^{\mathcal{N}_i} \Gamma_{ij} \left( S_{rf,i} \underline{\mathcal{T}}_{f,ij} - S_{rf,j} \underline{\mathcal{T}}_{f,ji} \right) \mathcal{V}_j \mathcal{V}_i, \tag{98}$$

$$\mathcal{L}_i = \mathcal{I}_i^s \mathcal{V}_i, \tag{99}$$

$$\mathcal{M}_i = \sum_{j=1}^{\mathcal{N}_i} (\underline{\mathcal{M}}_{ij} - \underline{\mathcal{M}}_{ji}) \mathcal{V}_j \mathcal{V}_i \tag{100}$$

$$\overline{\mathcal{M}}_i = \frac{1}{2} \sum_{j=1}^{\mathcal{N}_i} \left[ \underline{\mathcal{Y}}_{ij} \times \left( \underline{\mathcal{T}}_{ij} - \underline{\mathcal{T}}_{ji} \right) \right] \mathcal{V}_j \mathcal{V}_i$$

$$- \frac{1}{2} \sum_{j=1}^{\mathcal{N}_i} \left\{ \underline{\mathcal{Y}}_{ij} \times \left[ (1 - \Gamma_{ij}) \left( S_{r,i} \underline{\mathcal{T}}_{w,ij} - S_{r,j} \underline{\mathcal{T}}_{w,ji} \right) \right] \right\} \mathcal{V}_j \mathcal{V}_i$$

$$- \frac{1}{2} \sum_{j=1}^{\mathcal{N}_i} \left[ \underline{\mathcal{Y}}_{ij} \times \Gamma_{ij} \left( S_{rf,i} \underline{\mathcal{T}}_{f,ij} - S_{rf,j} \underline{\mathcal{T}}_{f,ji} \right) \right] \mathcal{V}_j \mathcal{V}_i, \tag{101}$$

where $\mathcal{A}$ is a linear global assembly operator [71], Vi and Vj are the volumes of material points i and j, respectively, and $\Gamma_{ij}$ is the indicator of the fracture point.

$$\Gamma_{ij} = \begin{cases} 1 & \text{if } i \text{ and } j \text{ are fracture points} \\ 0 & \text{otherwise} \end{cases} \tag{102}$$

In (97), (98), and (101), the effective force states and the fluid force states are written as

$$\underline{\mathcal{T}}_{ij} = \omega_{ij} \left( \overline{\boldsymbol{\sigma}}_i \mathcal{K}_i^{-1} \underline{\boldsymbol{\xi}}_{ij} + \underline{\alpha}_{1,ij} \mathcal{R}_{1,ij} \right), \tag{103}$$

$$\underline{\mathcal{T}}_{ji} = \omega_{ji} \left( \overline{\boldsymbol{\sigma}}_j \mathcal{K}_j^{-1} \underline{\boldsymbol{\xi}}_{ji} + \underline{\alpha}_{1,ji} \mathcal{R}_{1,ji} \right), \tag{104}$$

$$\underline{\mathcal{T}}_{w,ij} = \omega_{ij} p_{w,i} \mathbf{1} \mathcal{K}_i^{-1} \underline{\boldsymbol{\xi}}_{ij}, \tag{105}$$

$$\underline{\mathcal{T}}_{w,ji} = \omega_{ji} p_{w,j} \mathbf{1} \mathcal{K}_j^{-1} \underline{\boldsymbol{\xi}}_{ji}, \tag{106}$$

$$\underline{\mathcal{T}}_{f,ij} = \omega_{ij} p_{f,i} \mathbf{1} \mathcal{K}_i^{-1} \underline{\boldsymbol{\xi}}_{ij}, \tag{107}$$

$$\underline{\mathcal{T}}_{f,ji} = \omega_{ji} p_{f,j} \mathbf{1} \mathcal{K}_j^{-1} \underline{\boldsymbol{\xi}}_{ji}. \tag{108}$$

In (100), the moment states are written as

$$\underline{\mathcal{M}}_{ij} = \omega_{ij} \left( m_i \mathcal{K}_i^{-1} \underline{\boldsymbol{\xi}}_{ij} + \underline{\alpha}_{2,ij} \mathcal{R}_{2,ij} \right), \tag{109}$$

$$\underline{\mathcal{M}}_{ji} = \omega_{ji} \left( m_j \mathcal{K}_j^{-1} \underline{\boldsymbol{\xi}}_{ji} + \underline{\alpha}_{2,ji} \mathcal{R}_{2,ji} \right). \tag{110}$$

Similarly, the discretized mass balance equations can be written as

$$0 = \mathcal{A}_{i=1}^{\mathcal{P}} \left( \mathcal{X}_i + \mathcal{Q}_i + \widetilde{\mathcal{V}}_i + \mathcal{Q}_{s,i} \right), \tag{111}$$

$$0 = \mathcal{A}_{i=1}^{\mathcal{P}_f} \left( \mathcal{X}_{f,i} + \mathcal{Q}_{f,i} - \mathcal{Q}_{s,i} \right). \tag{112}$$

where Pf is the number of fracture points,

$$\mathcal{X}_i = -\phi_i \frac{\partial S_{r,i}}{\partial p_{w,i}} \dot{p}_{w,i} \mathcal{V}_i, \tag{113}$$

$$\mathcal{Q}_i = \frac{1}{\rho_w} \sum_{j=1}^{\mathcal{N}_i} (\underline{\mathcal{Q}}_{ij} - \underline{\mathcal{Q}}_{ji}) \mathcal{V}_j \mathcal{V}_i, \tag{114}$$

$$\widetilde{\mathcal{V}}_i = S_{r,i} \sum_{j=1}^{\mathcal{N}_i} \underline{\mathcal{V}}_i \mathcal{V}_j \mathcal{V}_i, \tag{115}$$

$$\mathcal{Q}_{s,i} = \mathcal{A}_i \left[ -\frac{k^r k_w}{\mu_w} \left( \frac{p_{w,i} - p_{f,i}}{l_{x,i}} \right) \right] / \mathcal{V}_i, \tag{116}$$

$$\mathcal{X}_{f,i} = -\phi_i \frac{\partial S_{r,fi}}{\partial p_{f,i}} \dot{p}_{f,i} \mathcal{V}_i, \tag{117}$$

$$\mathcal{Q}_{f,i} = \frac{1}{\rho_w} \sum_{j=1}^{\mathcal{N}_i} \Gamma_{ij} (\underline{\mathcal{Q}}_{f,ij} - \underline{\mathcal{Q}}_{f,ji}) \mathcal{V}_j \mathcal{V}_i. \tag{118}$$

In (114) and (118), the fluid flow states are written as

$$\underline{\mathcal{Q}}_{ij} = \underline{\omega}_{ij} \left( \rho_w \mathbf{q}_{w,i} \mathcal{K}_i^{-1} \underline{\xi}_{ij} \right) + \underline{\alpha}_{3,ij} \mathcal{R}_{w,ij}, \tag{119}$$

$$\underline{\mathcal{Q}}_{ji} = \underline{\omega}_{ji} \left( \rho_w \mathbf{q}_{w,j} \mathcal{K}_i^{-1} \underline{\xi}_{ji} \right) + \underline{\alpha}_{w,ji} \mathcal{R}_{3,ji}, \tag{120}$$

$$\underline{\mathcal{Q}}_{f,ij} = \underline{\omega}_{ij} \left( \rho_w \mathbf{q}_{f,i} \mathcal{K}_i^{-1} \underline{\xi}_{ij} \right) + \underline{\alpha}_{f,ij} \mathcal{R}_{f,ij}, \tag{121}$$

$$\underline{\mathcal{Q}}_{f,ji} = \underline{\omega}_{ji} \left( \rho_w \mathbf{q}_{f,j} \mathcal{K}_j^{-1} \underline{\xi}_{ji} \right) + \underline{\alpha}_{f,ij} \mathcal{R}_{f,ij}. \tag{122}$$

### 3.1.1. J integral

In this part, we compute the J-integral [47, 48] in the µPPM paradigm for a mode-I crack in two dimensions. Figure 4 presents a schematic of the discrete processing zone of a mode-I crack along the x direction and the path for the J-integral. It follows from (29) that the J-integral in the spatially discretized form can be written as

$$\begin{aligned}
\mathcal{J} = &\sum_{n=1}^{n_{\partial \mathcal{B}_1}} \mathcal{W} n_1 \Delta x - \sum_{i=1}^{n_{\mathcal{B}_1}} \sum_{j=1}^{n_{\mathcal{B}_2}} \left( \underline{\mathcal{T}}_{ij} \frac{\partial u_j}{\partial x_1} - \underline{\mathcal{T}}_{ji} \frac{\partial u_i}{\partial x_1} \right) \mathcal{V}_i \mathcal{V}_j / h \\
&- \frac{1}{2} \sum_{i=1}^{n_{\mathcal{B}_1}} \sum_{j=1}^{n_{\mathcal{B}_2}} \underline{\mathcal{X}}_{ij} \times \left( \underline{\mathcal{T}}_{ij} \frac{\partial \widehat{\omega}_j}{\partial x_1} - \underline{\mathcal{T}}_{ji} \frac{\partial \widehat{\omega}_i}{\partial x_1} \right) \mathcal{V}_i \mathcal{V}_j / h \\
&- \sum_{i=1}^{n_{\mathcal{B}_1}} \sum_{j=1}^{n_{\mathcal{B}_2}} \left( \underline{\mathcal{M}}_{ij} \frac{\partial \widehat{\omega}_j}{\partial x_1} - \underline{\mathcal{M}}_{ji} \frac{\partial \widehat{\omega}_i}{\partial x_1} \right) \mathcal{V}_i \mathcal{V}_j / h,
\end{aligned} \tag{123}$$

where n∂B1 is the number of material points on the boundary ∂B1, nB1 is the number of material points in the subregion B1 on a distance of 6 from the integral path, nB2 is the number of material points in the subregion B2 on a distance of 6 from the integral path, h is the thickness of crack, n is the unit normal to the path ∂B1, and x1 is the coordinate component parallel to the line crack (see Figure 4).

Similarly, from (32) the energy dissipation from the J-integral for the mode-I crack can be computed as

$$\mathcal{G} = \sum_{n=1}^{n_{\partial\mathcal{B}_1}} \mathscr{W} n_1 \Delta x - \sum_{i=1}^{n_{\mathcal{B}_1}} \sum_{j=1}^{n_{\mathcal{B}_2}} \left( \overline{\mathcal{T}}_{ij} \frac{\partial \boldsymbol{u}_j}{\partial x_1} - \overline{\mathcal{T}}_{ji} \frac{\partial \boldsymbol{u}_i}{\partial x_1} \right) \mathscr{V}_i \mathscr{V}_j / h$$
$$- \frac{1}{2} \sum_{i=1}^{n_{\mathcal{B}_1}} \sum_{j=1}^{n_{\mathcal{B}_2}} \mathscr{X}_{ij} \times \left( \overline{\mathcal{T}}_{ij} \frac{\partial \widehat{\omega}_j}{\partial x_1} - \overline{\mathcal{T}}_{ji} \frac{\partial \widehat{\omega}_i}{\partial x_1} \right) \mathscr{V}_i \mathscr{V}_j / h$$
$$- \sum_{i=1}^{n_{\mathcal{B}_1}} \sum_{j=1}^{n_{\mathcal{B}_2}} \left( \mathscr{M}_{ij} \frac{\partial \widehat{\omega}_j}{\partial x_1} - \mathscr{M}_{ji} \frac{\partial \widehat{\omega}_i}{\partial x_1} \right) \mathscr{V}_i \mathscr{V}_j / h. \tag{124}$$

In the section of numerical examples, (123) is used to compute the J-integral in the μPPM paradigm, which is compared with the closed-form solution for the J-integral in the classical micro-polar continuum model, and (30) is used to compute the mode-I crack propagation.

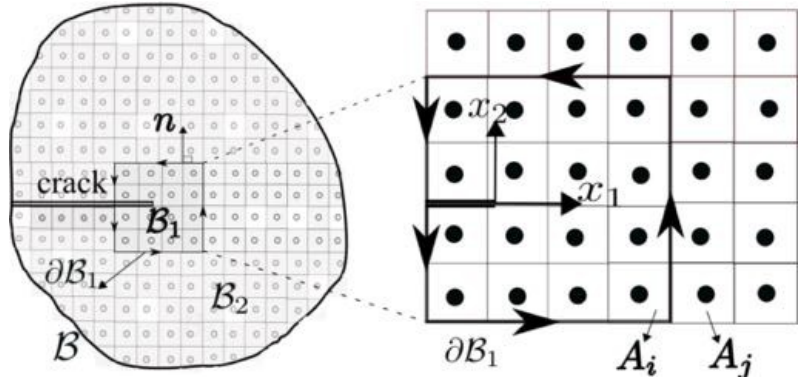

Figure 4: Schematic of the integration of J-integral of the μPPM paradigm for a mode-I crack.

### 3.2. Discretization in time

In this study, the fractional step/staggered explicit Newmark scheme [1, 71–73] is adopted to discretize the coupled μPPM paradigm in time. The fractional-step algorithm splits the coupled μPPM problem into a parallel deformation/fracturing model and an unsaturated/fracturing fluid flow problem. The energy balance check [73] is utilized to ensure the numerical stability of this algorithm for each solver in time. In what follows, we present the double-way staggered algorithm by introducing the solid deformation/fracturing solver, followed by the unsaturated/fracturing fluid flow solver. We note that, as stated at the beginning of this section, either the solid deformation/fracturing solver or the unsaturated/fracturing fluid flow solver can be called first, given the physical problem at hand.

#### 3.2.1. Solid deformation/fracturing solver under undrained conditions

Let $u_n$, $\dot{u}_n$ and $\ddot{u}_n$ be the displacement, velocity, and acceleration vectors at time step n and $\omega_n$, $\dot{\omega}_n$ and $\ddot{\omega}_n$ the mircorotation, angular velocity, and angular acceleration vectors at time step n. The predictors of displacement, micro-rotation, velocity, and angular velocity in a general Newmark scheme read

$$\tilde{\dot{u}}_{n+1} = \dot{u}_n + (1 - \beta_1) \Delta t \ddot{u}_n, \tag{125}$$
$$\tilde{\dot{\widehat{\omega}}}_{n+1} = \dot{\widehat{\omega}}_n + (1 - \beta_1) \Delta t \ddot{\widehat{\omega}}_n, \tag{126}$$
$$\tilde{u}_{n+1} = u_n + \Delta t \dot{u}_n + (1 - 2\beta_2) \frac{\Delta t^2}{2} \ddot{u}_n, \tag{127}$$
$$\tilde{\widehat{\omega}}_{n+1} = \widehat{\omega}_n + \Delta t \dot{\widehat{\omega}}_n + (1 - 2\beta_2) \frac{\Delta t^2}{2} \ddot{\widehat{\omega}}_n, \tag{128}$$

where $\beta_1$ and $\beta_2$ are numerical integration parameters. Given (125), (127), (126) and (128), T n+1, Tw, n+1, Mn+1 and M n+1 are determined from (127) and (128) from the classical local constitutive models introduced in Section through the µPPM correspondence principle. Then, the acceleration and angular acceleration at n + 1 are determined by

$$\ddot{u}_{n+1} = \mathcal{M}_{n+1}^{-1}\left(\widetilde{\mathscr{T}}_{n+1} - \widetilde{\mathscr{T}}_{w,n+1} - \mathcal{M}_{n+1}g\right), \tag{129}$$

$$\ddot{\omega}_{n+1} = \mathscr{L}_{n+1}^{-1}\left(\widetilde{\mathscr{M}}_{n+1} + \widetilde{\overline{\mathscr{M}}}_{n+1} + l_{n+1}V\right), \tag{130}$$

From (129) and (130), the velocity, angular velocity, displacement, and micro-rotation at time step n + 1 can be obtained as

$$\dot{u}_{n+1} = \dot{\widetilde{u}}_{n+1} + \beta_1 \Delta t \ddot{u}_{n+1}, \tag{131}$$

$$\dot{\hat{\omega}}_{n+1} = \dot{\widetilde{\omega}}_{n+1} + \beta_1 \Delta t \ddot{\hat{\omega}}_{n+1}, \tag{132}$$

$$u_{n+1} = \widetilde{u}_{n+1} + \beta_2 \Delta t^2 \ddot{u}_{n+1}, \tag{133}$$

$$\hat{\omega}_{n+1} = \widetilde{\hat{\omega}}_{n+1} + \beta_2 \Delta t^2 \ddot{\hat{\omega}}_{n+1}. \tag{134}$$

In this study, we adopt the explicit central difference solution scheme [71], i.e., $\beta_1$ = 1/2 and $\beta_2$ = 0. We note that the explicit method is efficient and robust to model dynamic problems [71]. The energy balance check is used to ensure numerical stability of the solid deformation/fracturing solver in time [73]. The internal energy, external energy, and kinetic energy of the solid deformation/fracturing process at time step n + 1 are defined as

$$\widehat{\mathscr{W}}_{int,n+1} = \widehat{\mathscr{W}}_{int,n} + \frac{\Delta t}{2}\left(\dot{u}_n + \frac{\Delta t}{2}\ddot{u}_n\right)\left[(\mathscr{T}_n - \mathscr{T}_{w_n}) + (\mathscr{T}_{n+1} - \mathscr{T}_{w_{n+1}})\right]$$
$$+ \frac{\Delta t}{2}\left(\dot{\hat{\omega}}_n + \frac{\Delta t}{2}\ddot{\hat{\omega}}_n\right)\left[(\mathscr{M}_n + \overline{\mathscr{M}}_n) + (\mathscr{M}_{n+1} + \overline{\mathscr{M}}_{n+1})\right], \tag{135}$$

$$\widehat{\mathscr{W}}_{ext,n+1} = \widehat{\mathscr{W}}_{ext,n} + \frac{\Delta t}{2}\left(\dot{u}_n + \frac{\Delta t}{2}\ddot{u}_n\right)(\mathcal{M}_n g + \mathcal{M}_{n+1}g)$$
$$+ \frac{\Delta t}{2}\left(\dot{\hat{\omega}}_n + \frac{\Delta t}{2}\ddot{\hat{\omega}}_n\right)(l_n \mathscr{V} + l_{n+1}\mathscr{V}), \tag{136}$$

$$\widehat{\mathscr{W}}_{kin,n+1} = \frac{1}{2}\dot{u}_{n+1}\mathcal{M}_{n+1}\dot{u}_{n+1} + \frac{1}{2}\dot{\hat{\omega}}_{n+1}\mathscr{L}_{n+1}\dot{\hat{\omega}}_{n+1}. \tag{137}$$

Then, it follows from the energy conservation criterion that

$$\left|\widehat{\mathscr{W}}_{kin,n+1} - \widehat{\mathscr{W}}_{ext,n+1} + \widehat{\mathscr{W}}_{int,n+1}\right| \leq \hat{\varepsilon}_1 \max\left(\widehat{\mathscr{W}}_{kin,n+1}, \widehat{\mathscr{W}}_{int,n+1}, \widehat{\mathscr{W}}_{ext,n+1}\right), \tag{138}$$

where $\hat{\varepsilon}_1$ is a small tolerance on the order of $10^{-2}$ [31, 73].

3.2.2. Unsaturated/fracturing fluid flow solver in the updated deformed configuration

The fluid pressure predictors at time step n + 1 through the general Newmark scheme [71] can be written as

$$\widetilde{p}_{w,n+1} = p_{w,n} + (1 - \beta_3)\Delta t \dot{p}_{w,n}, \tag{139}$$

$$\widetilde{p}_{f,n+1} = p_{f,n} + (1 - \beta_3)\Delta t \dot{p}_{f,n}. \tag{140}$$

Given (139) and (140), the fluid flow states can be determined. Then, p˙w,n+1 and p˙f,n+1 can be computed by

$$\dot{p}_{w,n+1} = -\left(\phi\frac{\partial\widetilde{s}_{r,n+1}}{\partial\widetilde{p}_{w,n+1}}\right)^{-1}\left(\widetilde{\mathcal{Q}}_{n+1} + \mathcal{V}_{n+1} + \widetilde{\mathcal{Q}}_{s,n+1}\right), \tag{141}$$

$$\dot{p}_{f,n+1} = \left(\phi\frac{\partial\widetilde{s}_{rf,n+1}}{\partial\widetilde{p}_{f,n+1}}\right)^{-1}\left(\widetilde{\mathcal{X}}_{f,n+1} + \widetilde{\mathcal{Q}}_{f,n+1} - \widetilde{\mathcal{Q}}_{s,n+1}\right), \tag{142}$$

It follows from (139), (140), (141), and (142), and the general Newmark scheme [71] that pw,n+1 and pf,n+1 can be written as

In this study, the explicit Newmark scheme, i.e., $\beta_3 = 1$, is adopted for the fluid flow solver. Similarly to the deformation/fracturing solver, the energy convergence criterion is used to ensure the numerical stability for the fluid flow solver. The internal energy, external energy, and kinetic energy of the fluid flow process at time step n + 1 are defined as

$$\widehat{\mathcal{W}}_{\text{int},n+1} = \widehat{\mathcal{W}}_{\text{int},n} + \frac{\Delta t}{2}\dot{p}_{w,n+1}\left[(\mathcal{Q}_n + \mathcal{V}_n) + (\mathcal{Q}_{n+1} + \mathcal{V}_{n+1})\right]$$
$$+ \frac{\Delta t}{2}\dot{p}_{f,n+1}\left[(\mathcal{X}_{f,n} + \mathcal{Q}_{f,n}) + (\mathcal{X}_{f,n+1} + \mathcal{Q}_{f,n+1})\right], \tag{145}$$

$$\widehat{\mathcal{W}}_{\text{ext},n+1} = \widehat{\mathcal{W}}_{\text{ext},n} + \frac{\Delta t}{2}\dot{p}_{w,n+1}(\mathcal{Q}_{s,n} + \mathcal{Q}_{s,n+1}) - \frac{\Delta t}{2}\dot{p}_{f,n+1}(\mathcal{Q}_{s,n} + \mathcal{Q}_{s,n+1}). \tag{146}$$

Then it follows from the energy conservation criterion that

$$|\widehat{\mathcal{W}}_{\text{int},n+1} - \widehat{\mathcal{W}}_{\text{ext},n+1}| \leq \hat{\varepsilon}_2 \max\left(\widehat{\mathcal{W}}_{\text{int},n+1}, \widehat{\mathcal{W}}_{\text{ext},n+1}\right), \tag{147}$$

where $\hat{\varepsilon}_2$ is a small tolerance on the order of $10^{-2}$ [31, 73]. It is worth noting that following the lines in [29], a fractional-step implicit-implicit algorithm can be formulated to solve the coupled μPPM paradigm in time, which is beyond the scope of the present study. For the numerical implementation algorithm for micro-polar plasticity models, we refer to the celebrated literature on the subject (e.g., [42, 44, 69, 74, 75]).

## 4. Numerical examples

In this section, we present four numerical examples to evaluate the implemented μPPM paradigm and to demonstrate its efficacy and capability of modeling shear banding and cracking in porous media. Example 1 deals with the mode I crack driven by fluid pressure or solid deformation. Example 2 validates and verifies the proposed μPPM J integral for the mode I and II cracks. Example 3 deals with the shear banding of unsaturated porous media under the biaxial compression test. Example 4 deals with the propagation of the mode-I crack under displacement control loading in unsaturated porous media.

*4.1. Example 1: Mode I crack*

This example simulates the mode I crack driven by either fluid pressure or solid deformation. We first present the case driven by fluid pressure, followed by the case driven by solid deformation.

**Algorithm 1** Fractional-step explicit-explicit algorithm for the $\mu$PPM paradigm

1: Given: $u_n, \dot{u}_n, \ddot{u}_n, \widehat{\omega}_n, \dot{\widehat{\omega}}_n, \ddot{\widehat{\omega}}_n, p_{w,n}, \dot{p}_{w,n}, p_{f,n}, \dot{p}_{f,n}, t_n, \Delta t$
2: Compute: $u_{n+1}, \dot{u}_{n+1}, \ddot{u}_{n+1}, \widehat{\omega}_{n+1}, \dot{\widehat{\omega}}_{n+1}, \ddot{\widehat{\omega}}_{n+1}, p_{w,n+1}, \dot{p}_{w,n+1}, p_{f,n+1}, \dot{p}_{f,n+1}$
3: Update time $t_{n+1} = t_n + \Delta t$
4: **while** $t_{n+1} \leq t_{final}$ **do**
5:     Compute the velocity predictor $\tilde{\dot{u}}_{n+1}$ using (125)
6:     Compute the micro rotation rate predictor $\tilde{\dot{\widehat{\omega}}}_{n+1}$ using (126)
7:     Apply boundary conditions
8:     Compute displacement predictor $\tilde{u}_{n+1}$ using (127)
9:     Compute micro rotation predictor $\tilde{\widehat{\omega}}_{n+1}$ using (128)
10:     Compute effective force, fluid force and moment vectors via Algorithm 2
11:     Compute $\widetilde{\mathcal{M}}_{n+1}$ and $\widetilde{\mathcal{L}}_{n+1}$
12:     Solve the accelerations $\ddot{u}_{n+1}$ using (129) and $\ddot{\widehat{\omega}}_{n+1}$ using (130)
13:     Update velocity $\dot{u}_{n+1}$ using (131) and micro rotation rate $\dot{\widehat{\omega}}_{n+1}$ using (132)
14:     Update displacement $u_{n+1}$ using (133) and micro rotation $\widehat{\omega}_{n+1}$ using (134)
15:     Compute kinematic energy $\widehat{\mathcal{W}}_{kin,n+1}$ using (137)
16:     Compute internal energy $\widehat{\mathcal{W}}_{int,n+1}$ using (135) and external energy $\widehat{\mathcal{W}}_{ext,n+1}$ using (136)
17:     Check energy balance

18:     Compute J integral for J integral criterion
19:     Update the list of broken bonds fracture points via Algorithm 3
20:     Compute bond energy for energy based bond breakage criterion
21:     Update the list of broken bonds and fracture points via Algorithm 4

22:     Compute water pressure predictor $\tilde{p}_{w,n+1}$ using (139)
23:     Compute fracture pressure predictor $\tilde{p}_{f,n+1}$ using (140)
24:     Compute the fluid flow in the bulk and fracture points via Algorithm 5
25:     Solve water pressure rate $\dot{p}_{w,n+1}$ using (141)
26:     Solve fracture pressure rate $\dot{p}_{f,n+1}$ using (142)
27:     Update water pressure $p_{w,n+1}$ using (143) and fracture pressure $p_{f,n+1}$ using (144)
28:     Compute internal energy $\widehat{\widehat{\mathcal{W}}}_{int,n+1}$ using (145) and external energy $\widehat{\widehat{\mathcal{W}}}_{ext,n+1}$ using (146)
29:     Check energy balance
30: **end while**
31: $n \leftarrow n + 1$

---

**Algorithm 2** Compute effective force, fluid force and moment states

1: Given: displacement predictor $\tilde{u}_{n+1}$, micro rotation predictor $\tilde{\widehat{\omega}}_{n+1}$, water pressure predictor at bulk point $\tilde{p}_{w,n+1}$ and at fracture point $\tilde{p}_{f,n+1}$, $\underline{\omega}_n$, $\Gamma_n$
2: Compute: effective force vector, moment vector and fluid force vector
3: **for** all points **do**
4:     Compute composite state for solid phase
5:     Compute relative rotation state for solid phase
6:     Compute shape tensor $\boldsymbol{K}$
7:     Compute effective stress $\bar{\boldsymbol{\sigma}}$ and couple stress $\boldsymbol{m}$
8:     Compute force states $\underline{\bar{\mathcal{T}}}$ using (103), $\underline{\mathcal{T}}_w$ using (105) and $\underline{\mathcal{T}}_f$ using (107)
9:     Compute moment state $\underline{\mathcal{M}}$ using (109)
10:     Compute effective force vector $\widetilde{\bar{\mathcal{T}}}_{n+1}$ using (97)
11:     Compute fluid force vector $\widetilde{\mathcal{T}}_{w,n+1}$ using (98)
12:     Compute moment vector $\widetilde{\mathcal{M}}_{n+1}$ using (100) and $\widetilde{\mathcal{M}}_{\underline{\mathcal{T}}_{n+1}}$ using (101)
13: **end for**

**Algorithm 3** Compute J Integral and update the fracture point

1: Given: $\overline{\mathcal{T}}_{n+1}$, $u_{n+1}$, $\mathcal{M}_{n+1}$ and $\widehat{\omega}_{n+1}$
2: Compute J integral and update influence function $\underline{\omega}$ and update fracture points list
3: Define crack path direction $x_1$
4: Define the integral contour path
5: **for** all points along the contour **do**
6:     Compute strain energy density
7:     Compute first term of J integral at (123)
8: **end for**
9: **for** all points in regions $\mathcal{B}_1$ and $\mathcal{B}_2$ **do**
10:     Compute the second, third and 4th terms of J integral at (123)
11: **end for**
12: Compute J integral
13: **if** $J > J_{cr}$ **then**
14:     **for** all points **do**
15:         **for** each neighbour **do**
16:             **if** the bond cross the new crack surface **then**
17:                 Update influence function $\underline{\omega}$
18:             **end if**
19:         **end for**
20:         Update damage variable $D$
21:         **if** $D > D_{cr}$ **then**
22:             Set material point as fracture point
23:         **end if**
24:     **end for**
25: **end if**

**Algorithm 4** Compute bond energy and update the fracture points

1: Given: $\overline{\mathcal{T}}_{n+1}$, $\Delta\mathcal{U}$, $\mathcal{M}_{n+1}$ and $\Delta\Omega$
2: Compute bond energy and update influence function $\underline{\omega}$ and update fracture points list
3: **for** all points **do**
4:     **for** each neighbor **do**
5:         Compute bond energy $\mathcal{W}$
6:         **if** $\mathcal{W} > \mathcal{W}_{cr}$ **then**
7:             Update influence function $\underline{\omega}$
8:             Update damage variable $D$
9:             Sum the energy in the bond to total energy dissipated at the point
10:         **end if**
11:     **end for**
12:     **if** $D > D_{cr}$ **then**
13:         Set material point as fracture point
14:     **end if**
15: **end for**

**Algorithm 5** Compute fluid flow in the bulk and fracture
1: Given: $\dot{\boldsymbol{u}}_{n+1}\ \widetilde{p}_{w,n+1},\ \widetilde{p}_{f,n+1},\ \Gamma_{n+1}$
2: Compute: bulk fluid flow and fracture fluid flow vectors
3: **for** all points **do**
4:     Compute pressure potential state
5:     Compute shape tensor $\boldsymbol{K}$
6:     Compute pressure gradient state $\widetilde{\nabla\Phi}$
7:     Compute relative permeability $k^r$
8:     Compute flux vector $\boldsymbol{q}_w$
9:     Compute flow state $\underline{\mathcal{Q}}$
10:    Compute bulk fluid flow $\widetilde{\mathcal{Q}}_{n+1}$ using (114)
11: **end for**
12: **for** all points **do**
13:    Compute fracture width $a_f$
14:    Compute fracture permeability $k_f$
15:    **for** each neighbor **do**
16:       Compute fracture pressure states
17:       Compute fracture pressure gradient state $\widetilde{\nabla\Phi}_f$
18:       Compute fracture flow state $\underline{\mathcal{Q}}_f$
19:    **end for**
20:    Compute fracture fluid flow $\widetilde{\mathcal{Q}}_{f,n+1}$ using (118)
21:    Compute source term $\widetilde{\mathcal{Q}}_{s,n+1}$ using (116)
22: **end for**

*4.1.1. Case 1: Fluid driven crack*
*4.1.2.*

This example simulates the fluid pressure-driven crack in a saturated elastic porous material. Figure 6 plots the model setup for case 1. The initial crack length is 0.05 m. The fluid flow rate $q = 5 \times 10^{-5}$ m2/s is imposed on the left end of the crack. The water pressure is prescribed zero on the top, bottom, and right boundaries of the specimen. The horizontal displacement is fixed on the right boundary. The displacement in the vertical direction is fixed on the top and bottom boundaries. The elastic micro-polar constitutive model for the solid introduced in Section 2.4.1 is adopted. Following [76], the input material parameters are: bulk modulus K = 270.83 MPa, shear modulus μ = 220.33 MPa, micropolar shear modulus μc = 100 MPa, micropolar length scale $l = 7.5 \times 10^{-3}$ m, solid density $\rho_s$= 2000 kg/m3, initial porosity $\phi_0$= 0.2, water density $\rho_w$= 1000 kg/m3, and hydraulic conductivity kw = 2.78 $\times 10^{-9}$ m/s. The water pressure is assumed incompressible. The stabilization parameter G = 0.5 is used. For this case, $G_{cr}$= 95 N/m is used for the energy-based bond breakage criterion [76]. The loading time is t = 0.5 s, and the time increment is $\Delta t = 5 \times 10^{-6}$ s. In what follows, we present the numerical results.

Figure 7 plots the snapshots of the damage variable contour in the deformed configuration at three loading times. Figure 8 compares the contours of water pressure in the deformed con- figuration at three loading times with the results in [76]. The contours of water pressure and the damage variable are similar to the numerical results in [76](see Figures 7 and 8). Figure 9 presents the snapshots of micro rotation in the deformed configuration at the three loading times. As shown in Figure 9, the micro-rotation of material points is concentrated on the crack tip, and the magnitude of the micro-rotation of those material points increases as the crack propagates. Next, we present the results of the simulations with two spatial discretizations. For the two cases, the specimen is discretized into 66×66 uniform points ($\Delta x_1 = 3.75 \times 10^{-3}$ m) and $100 \times 100$ uniform points ($\Delta x_2 = 2.5 \times 10^{-3}$ m), respectively. The same micro-length scale and horizon, i.e., ($\delta = 7.5 \times 10^{-3}$ m, is assumed. The other conditions remain the same. Figure 10 plots the contours of the damage variable in the deformed configuration at loading ($\bar{\sigma}$= 1 MPa for both discretizations. Figure 11 shows the contours of the water pressure in the deformed configuration at loading ($\bar{\sigma}$ = 1 MPa for both discretizations. Figure 12 plots the contours of the micro rotation of material points at loading ($\bar{\sigma}$= 1 MPa for the two discretizations. The results in these figures have demonstrated that with the same horizon, the numerical results are less influenced by the spatial discretization scheme.

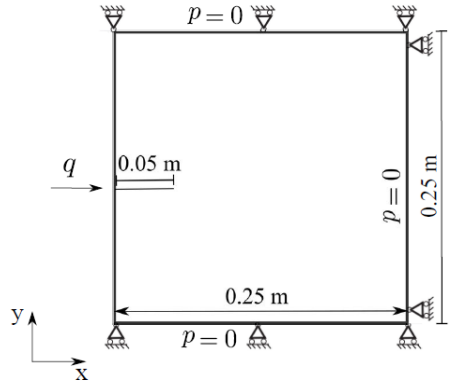

Figure 6: Problem setup for the case of fluid pressure driven crack.

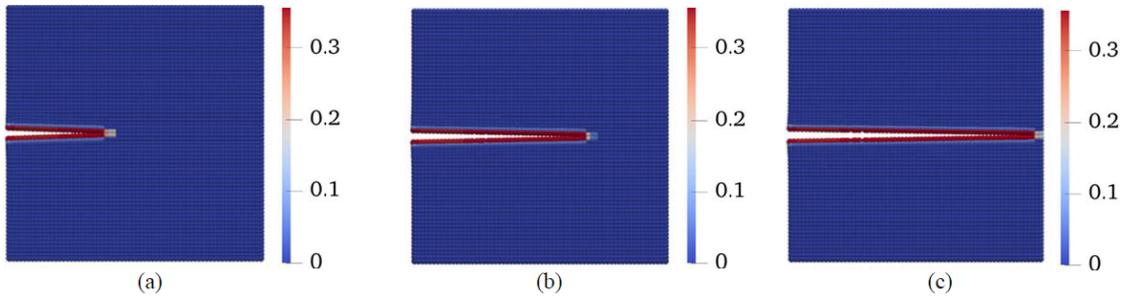

Figure 7: Contours of the damage variable in the deformed configuration (magnification factor (MF) =100) at (a) $t_1 = 0.3$ s, (b) $t_2 = 0.4$, and s (c) $t_3 = 0.5$ s.

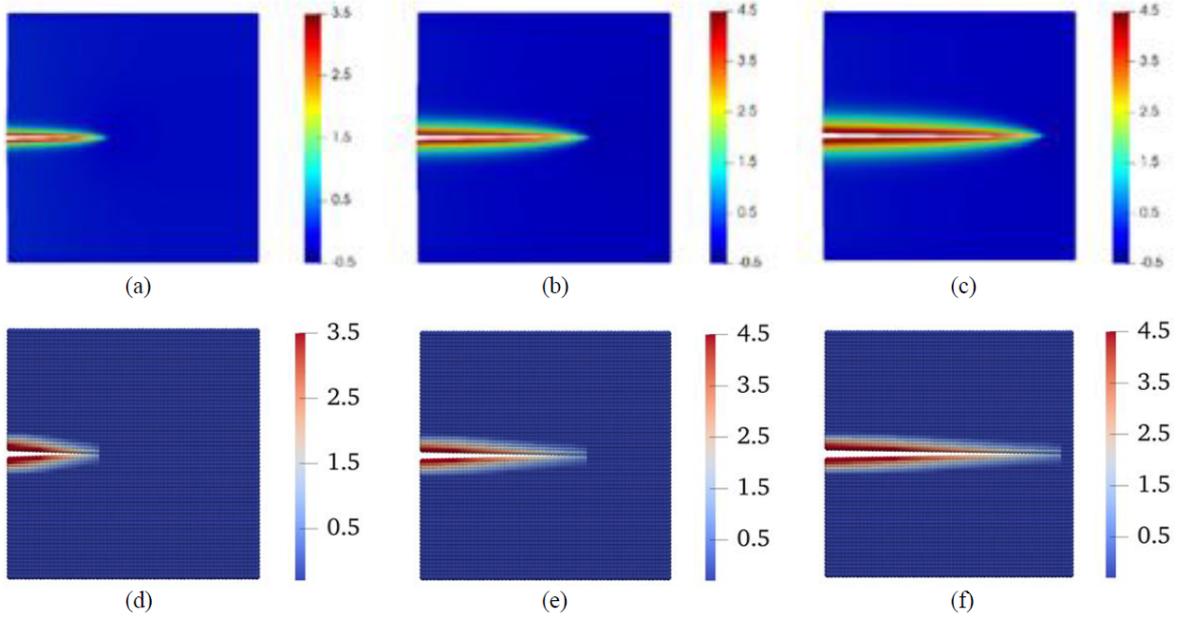

Figure 8: Contours of water pressure (MPa) in the deformed configuration (MF = 100) from [76] at (a) $t_1 = 0.3$ s, (b) $t_2 = 0.4$ s, and (c) $t_3 = 0.5$ s and from our numerical model at (d) $t_1 = 0.3$ s, (e) $t_2 = 0.4$ s, and (f) $t_3 = 0.5$ s.

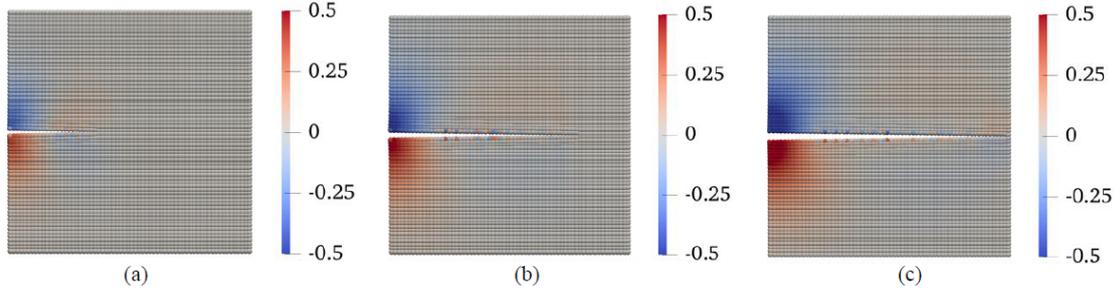

Figure 9: Contours of the micro-rotation (degree) in the deformed configuration (MF = 100) at (a) $t_1 = 0.3$ s, (b) $t_2 = 0.4$ s, and (c) $t_3 = 0.5$ s.

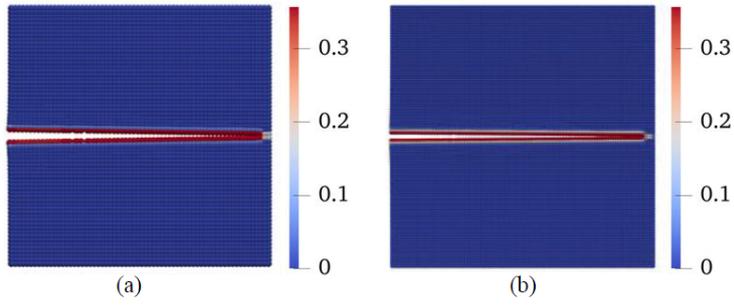

Figure 10: Contours of the damage variable on the deformed configuration (MF = 100) at loading $\bar{\sigma} = 1$ MPa for cases ($\delta = 7.5 \times 10^{-3}$ m): (a) $\Delta x_1 = 3.75 \times 10^{-3}$ m and (b) $\Delta x_2 = 2.5 \times 10^{-3}$ m.

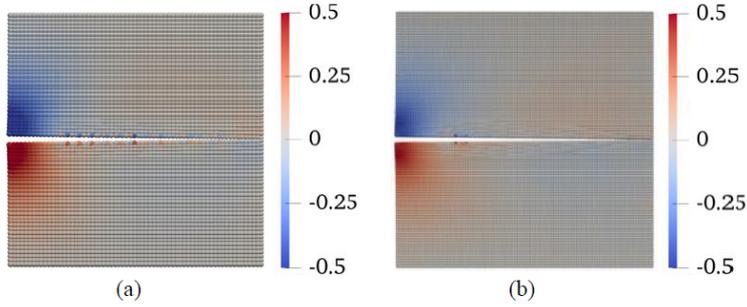

Figure 11: Contours of water pressure (MPa) in deformed configuration (MF = 100) at loading $\bar{\sigma} = 1$ MPa for cases ($\delta = 7.5 \times 10^{-3}$ m): (a) $\Delta x_1 = 3.75 \times 10^{-3}$ m and (b) $\Delta x_2 = 2.5 \times 10^{-3}$ m.

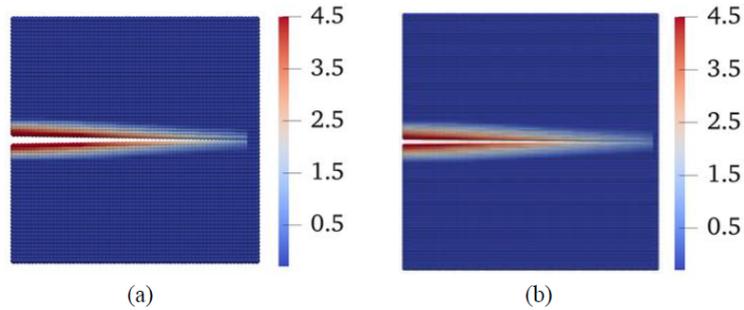

Figure 12: Contours of micro rotation (degree) in deformed configurations (MF = 100) at loading $\bar{\sigma} = 1$ MPa for cases ($\delta = 7.5 \times 10^{-3}$ m): (a) $\Delta x_1 = 3.75 \times 10^{-3}$ m and (b) $\Delta x_2 = 2.5 \times 10^{-3}$ m.

*4.1.3. Case 2: Deformation driven crack*

Case 2 deals with the deformation-driven crack in the saturated elastic porous material, as assumed in case 1. Figure 13 plots the dimensions, boundary, and loading conditions for case 2. We note that the dimensions of the specimen and the initial crack are the same as in case 1, except for the loading condition. All boundaries are impervious. The input material parameters are the same as in case 1. The vertical load is applied on the top and bottom boundaries through a ramp function from 0 to 1 MPa. The specimen is discretized into 100 × 100 uniform material points with $\Delta x = 2.5 \times 10^{-3}$ m. The horizon size is the same as the micro-polar length scale, i.e., $\delta = 3.06 \Delta x$. The simulation time is $t = 5 \times 10^{-3}$ s. The time increment is $\Delta t = 5 \times 10^{-7}$ s.

Figure 14 plots the snapshots of the damage variable at ($\bar{\sigma}$= 0.6 MPa, 0.8 MPa, and 1MPa, respectively. As shown in Figure 14, at the load ($\bar{\sigma}$= 0.6 MPa, the crack starts to grow, and at the load ($\bar{\sigma}$= 1 MPa, the crack length reaches 0.1975 m. Figure 15 shows the snapshots of water pressure at the same three loading stages. Figure 15 shows that the maximum water pressure is at the crack tip. It is implied from Figure 15 that the water pressure in the material points around the crack path is zero due to fluid flow to the fracture points. Figure 16 presents the snapshots of micro rotation at three loading stages of the simulation. Figure 9 shows that the micro rotation of material points at the crack tip is the maximum. It is implied from the results that the magnitude of the micro rotation of material points increases as the crack propagates.

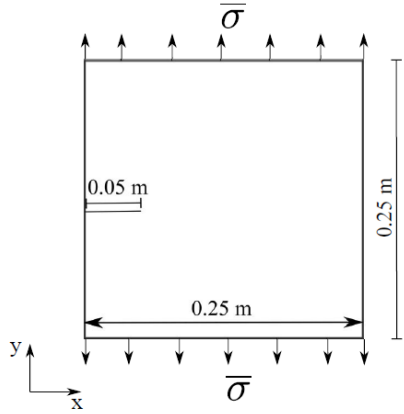

Figure 13: Problem setup for the deformation driven crack.

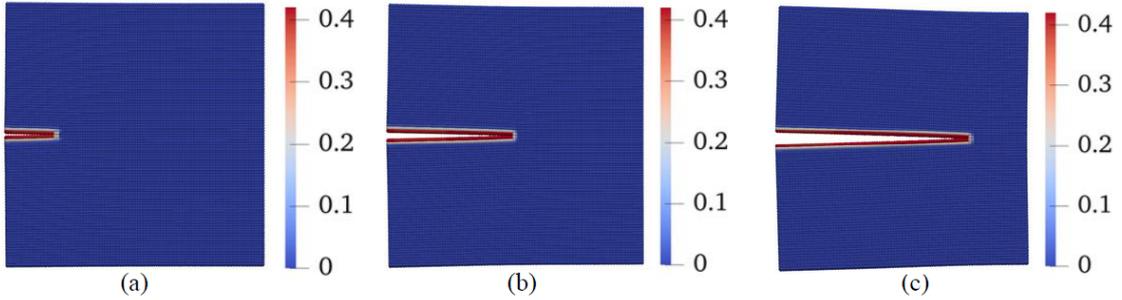

Figure 14: Contours of the damage variable in deformed configurations (MF = 100) at (a) $\bar{\sigma}_1 = 0.6$ MPa, (b) $\bar{\sigma}_2 = 0.8$ MPa, and (c) $\bar{\sigma}_3 = 1$ MPa.

Next, we present the results of the simulations with two spatial discretization assuming the same length scale. For the two cases, the specimen is discretized into 66×66 uniform points ($\Delta x_1 = 3.75 \times 10^{-3}$ m) and 100 × 100 uniform points ($\Delta x_1 = 2.5 \times 10^{-3}$ m), respectively. The same micro-length scale and horizon are assumed, i.e., $\delta = 7.5 \times 10^{-3}$ m. The other conditions are the same. Figure 17 plots the contours of the damage variable in the deformed configuration at loading ($\bar{\sigma}$ = 1 MPa for both discretization. Figure 18 shows the contours of the water pressure in the deformed configuration at loading ($\bar{\sigma}$ = 1 MPa for both discretization. Figure 19 demonstrates the contours of the micro rotation of material points in the deformed configuration at loading ($\bar{\sigma}$ = 1 MPa for both discretization. It can be implied from the results that with the same horizon, the numerical results are less influenced by spatial discretization.

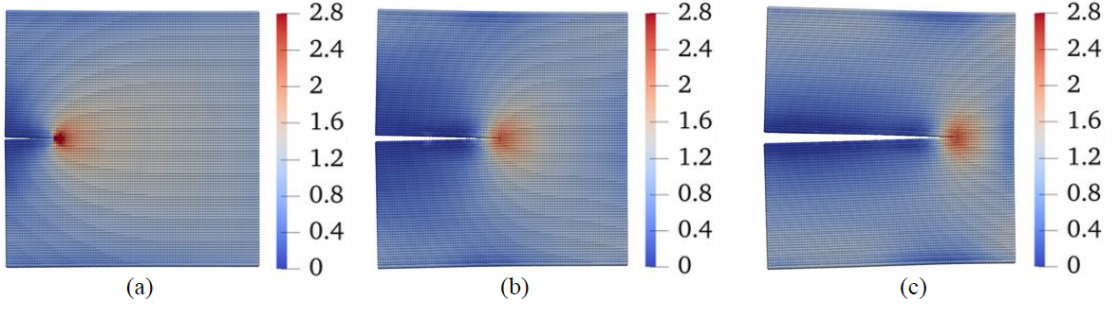

Figure 15: Contours of water pressure (MPa) in deformed configurations (MF = 100) at (a) $\overline{\sigma}_1 = 0.6$ MPa, (b) $\overline{\sigma}_2 = 0.8$ MPa, and (c) $\overline{\sigma}_3 = 1$ MPa.

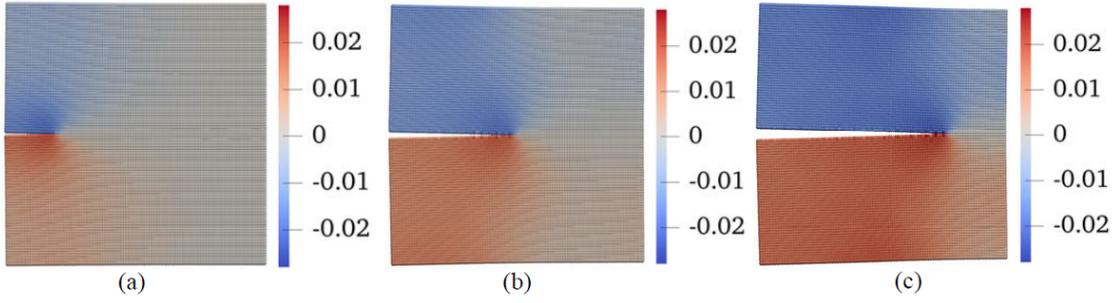

Figure 16: Contours of the micro rotation (degree) in deformed configurations (MF = 100) at loading (a) $\overline{\sigma}_1 = 0.6$ MPa, (b) $\overline{\sigma}_2 = 0.8$ MPa, and (c) $\overline{\sigma}_3 = 1$ MPa.

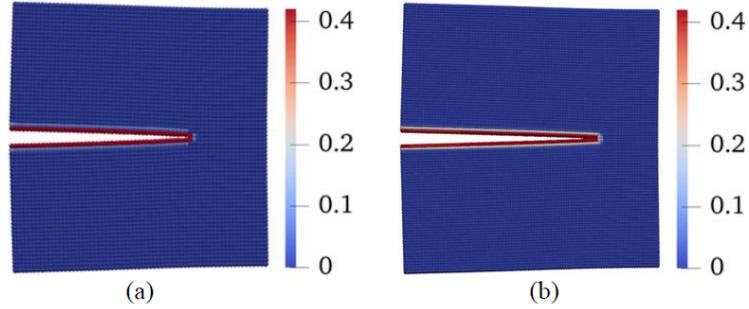

Figure 17: Contours of damage variable in deformed configuration (MF = 100) at loading $\overline{\sigma} = 1$ MPa for cases ($\delta = 7.5 \times 10^{-3}$ m): (a) $\Delta x_1 = 3.75 \times 10^{-3}$ m and (b) $\Delta x_2 = 2.5 \times 10^{-3}$ m.

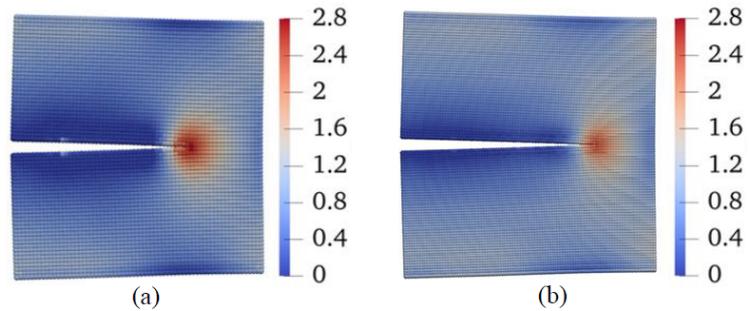

Figure 18: Contours of water pressure (MPa) in deformed configuration (MF = 100) at loading $\overline{\sigma} = 1$ MPa for cases ($\delta = 7.5 \times 10^{-3}$ m): (a) $\Delta x_1 = 3.75 \times 10^{-3}$ m and (b) $\Delta x_2 = 2.5 \times 10^{-3}$ m.

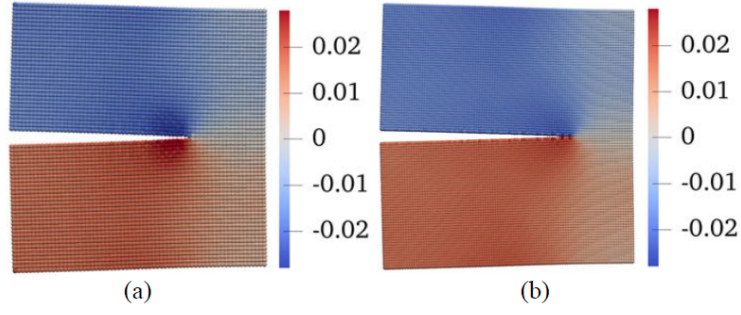

Figure 19: Contours of micro rotation (degree) in deformed configurations (MF = 100) at loading $\bar{\sigma} = 1$ MPa for cases ($\delta = 7.5 \times 10^{-3}$ m): (a) $\Delta x_1 = 3.75 \times 10^{-3}$ m and (b) $\Delta x_2 = 2.5 \times 10^{-3}$ m.

*4.2. Example 2: Validation of J-integral*

This example deals with the J integral computation in the µPPM framework for the mode I and II cracks in an elastic porous material. In this example, we compare our numerical results of J-integral with the classical micropolar solution [45]. Furthermore, we study the impact of the micropolar length scale, horizon, and spatial discretization schemes on the J integral computation.

*4.2.1. J-integral: Mode I crack*

In this part, we compute the J-integral of a mode I crack. Figure 20 plots the model setup. The initial crack length is 0.05 m. A tensile stress ($\bar{\sigma}$ = 1 MPa is applied on the top and bottom boundaries through a ramp function. The elastic material properties are: bulk modulus K = 60 GPa, shear modulus µ = 27.7 GPa, micropolar shear modulus µc = 14 GPa, solid density s = 3000 kg/m3, and initial porosity ¢0 = 0.3. The stabilization parameter G = 0.5 is used. The specimen is discretized into 100⇥100 material points with a uniform grid size $\Delta x = 1 \times 10^{-3}$ m. The simulation time is $t = 1 \times 10^{-2}$ s and the time increment $\Delta t = 2.5 \times 10^{-7}$ s.

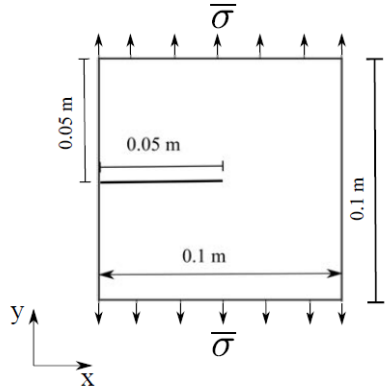

Figure 20: Model setup for the mode I crack.

We first show the results of a base simulation. Figure 21 plots the contour of the damage variable, the micro rotation, and the displacement norm superimposed on the deformed configuration at the final loading step. As shown in Figure 21, the magnitude of micro rotation at the crack tipis very small for the mode I crack.

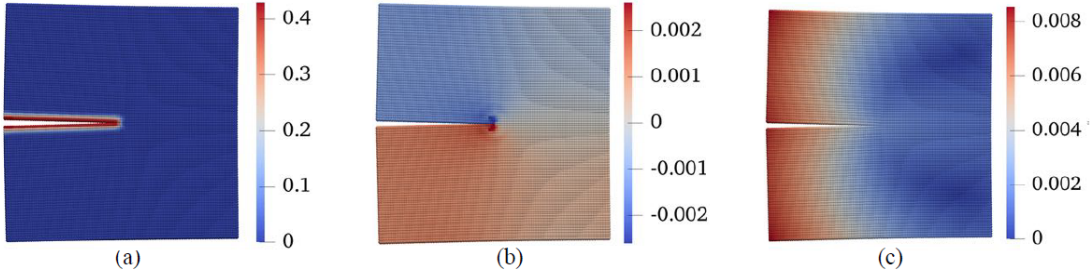

Figure 21: (a) Contour of the damage variable, (b) Contour of micro rotation (degree), and (c) Contour of displacement norm ($\times 10^{-3}$ m) in the deformed configuration (MF = 200) at the final load step.

The micropolar length scale is hypothesized to affect the J integral for the mode I crack. To test this hypothesis, we compute the J-integral with three micropolar length scales, $l = 3 \times 10^{-3}$ m, $4 \times 10^{-3}$ m, and $5 \times 10^{-3}$ m. The horizon size is the same as the micropolar length scale. Table 1 summarize the results from μPPM and the classical micropolar solution [45]. In Table 1, we also present the translational and micro-rotational parts of the J integral of our numerical results for each micro-polar length scale. The results in Table 1 show that our numerical solution is slightly greater than the closed-form solution of the classical micro-polar formulation. With the increase of the micro-polar length scale, the J-integral for the mode I crack decreases for both the μPPM and classical micro-polar solutions. Table 1 shows that the micro rotational part of the J integral is much smaller than the translational displacement part due to the mode I crack.

Table 1: Comparison of the J-integral from the μPPM solution and the classical micro-polar solution for the mode I crack assuming three micro-polar length scales.

| $l$ ($\times 10^{-3}$ m) | μPPM | Translational part | Micro-rotational part | Classical micropolar |
|---|---|---|---|---|
| | | J integral (Pa. m) | | |
| 3 | 21.64 | 21.73 | 0.21 | 19.81 |
| 4 | 20.48 | 20.3 | 0.18 | 19.7 |
| 5 | 20.04 | 19.89 | 0.15 | 19.66 |

Next, we study the impact of the grid size on the J-integral result. For this purpose, we present the results of J integral with four spatial discretization schemes, i.e., $\Delta x = 1.33 \times 10^{-3}$ m, $1 \times 10^{-3}$ m and $0.8 \times 10^{-3}$ m. The same horizon $\delta = 4 \times 10^{-3}$ m and the same internal length scale $l = 4 \times 10^{-3}$ m are adopted. Table 2 presents the results of the J-integral obtained from the four grid sizes. Table 2 shows that the grid size slightly impacts the J-integral. This impact may be decreased using a smaller grid size under the same micro-polar length. In the part, we present the results of J-integral for the mode II crack.

Table 2: Summary of the J-integrals from μPPM with 4 grid sizes for the mode I crack.

| $\Delta x$ ($\times 10^{-3}$ m) | $l$ ($\times 10^{-3}$ m) | J integral (Pa. m) |
|---|---|---|
| 1.33 | 4 | 21.24 |
| 1 | 4 | 20.48 |
| 0.8 | 4 | 20.01 |
| 0.67 | 4 | 19.93 |

*4.2.2. J-integral: Mode II crack*

We compute the μPPM J-integral for the mode II crack in this part. Figure 22 plots the model setup for the mode II crack. The initial crack length is 0.05 m. A shear force is applied on the top boundary through a ramp function from zero to 1 kN. The material parameters are the same for the mode I crack in this example. The specimen is discretized into $100 \times 100$ material points with a uniform grid size $\Delta x = 1 \times 10^{-3}$ m. The simulation time is $t = 1 \times 10^{-2}$ s and $\Delta t = 2.5 \times 10^{-7}$ s. We first present the base simulation. Figure 23 plots the contour of the damage variable, the micro rotation, and the displacement magnitude in the deformed configuration at the final load step. Figure 23 shows that the micro-rotation magnitude at the crack tip is more significant in the mode II crack than that for the mode I crack case. This could be due to the micro-rotation caused by the shear loading on the top boundary.

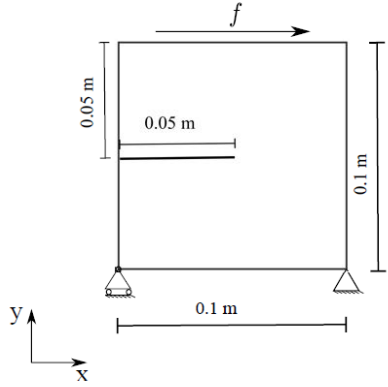

Figure 22: The problem setup for Mode II crack of example 2

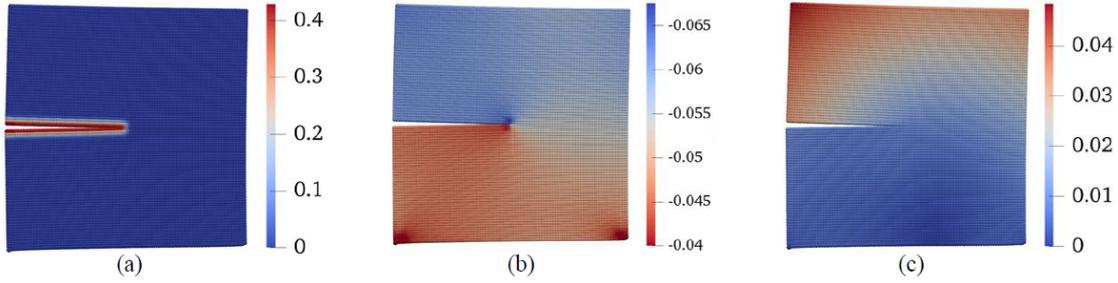

Figure 23: (a) Contour of the damage variable, (b) contour of the micro rotation (degree), and (c) contour of the displacement norm ($\times 10^{-3}$m) in the deformed configuration (MF = 100) at the last loading step.

Similarly, to the mode I crack, it is postulated that the micro-polar length affects the J-integral in the mode II crack. To test this postulation, we compute the J-integral of the mode II crack with three micro-polar length scales, $l = 3 \times 10^{-3}$m, $4 \times 10^{-3}$ m, and $5 \times 10^{-3}$ m. The horizon equals the micropolar length scale. Table 3 summarizes the results from µPPM and the classical micro-polar solution [45]. Table 3 also presents the translational and micro-rotational parts of the J integral of the mode II crack for different micro-polar length scales. The results show that our numerical solutions are close to the solution of the classical micropolar formulation. Moreover, the results in Table 3 also show that the micro-rotational part of the J integral for the mode II crack could play a more important role in the J-integral than the mode I crack.

Table 3: Comparison of J-integral from µPPM and the classical micropolar model for the mode II crack.

| $l$ ($\times 10^{-3}$m) | µPPM | Translational part | Micro-rotational part | Classical micro-polar |
|---|---|---|---|---|
| | | J integral (Pa. m) | | |
| 3 | 11.76 | 7.14 | 4.62 | 10.26 |
| 4 | 10.51 | 6.49 | 4.02 | 9.94 |
| 5 | 10.21 | 6.32 | 3.69 | 9.71 |

Next, we present the results of the J integral computed from four spatial discretization schemes, i.e., $\Delta x = 1.33 \times 10^{-3}$ m, $1 \times 10^{-3}$ m and $0.8 \times 10^{-3}$m with the same internal length scale and horizon, i.e., $4 \times 10^{-3}$ m. Table 2 summarizes the results. The results show that the J integral value converges with decreasing the grid size.

Table 4: Summary of the J-integrals from μPPM with 4 grid sizes for the mode II crack.

| $\Delta x$ ($\times 10^{-3}$ m) | $l$ ($\times 10^{-3}$ m) | J integral (Pa. m) |
|---|---|---|
| 1.33 | 4 | 11.21 |
| 1 | 4 | 10.51 |
| 0.8 | 4 | 10.12 |
| 0.67 | 4 | 10.05 |

*4.3. Example 3: Conjugate shear banding under bi-axial loading*

This example deals with the conjugate shear banding in unsaturated plastic porous media under dynamic loading conditions. Figure 24 plots the model setup for this example. A vertical displacement $u_y = 10 \times 10^{-3}$ m is applied on the top boundary with a rate $\dot{u}_y$ = 0.1 m/s. The constant lateral confining pressure of 0.1 MPa is applied on the left and right boundaries. All boundaries are impervious to fluids. The micro-plastic model in Section 2.4.2 is adopted for this example. For the base simulation, the input material parameters are as follows: solid phase density s = 2000 kg/m3, bulk modulus K = 27.8 MPa, shear modulus μ = 20.8 MPa, micropolar shear modulus μc = 40.6 MPa, micropolar length scale l = 45 × $10^{-3}$ m, initial porosity $\phi_0$ = 0.35, water viscosity μw = 1 × $10^{-3}$ Pa s, hydraulic conductivity kw = 1 × $10^{-9}$ m/s, n = 1.8, sa = 50 kPa, initial cohesion $c_0$ = 0.8 MPa, residual cohesion $c_r$ = 0.25 MPa, linear softening modulus h = -10 MPa, dilatation angle 35° and frictional angle 35°. The stabilization parameter is G = 0.1. The initial water pressure is $p_0$ = 25 kPa. The specimen is discretized into 60×120 material points through a uniform grid with $\Delta x = 1.67 \times 10^{-3}$ m. The micro-polar length scale and the horizon size are the same, i.e., $\delta$ = 3.056$\Delta x$. The simulation time is t = 0.1 s and the time increment is $\Delta t = 2.5 \times 10^{-6}$ s.

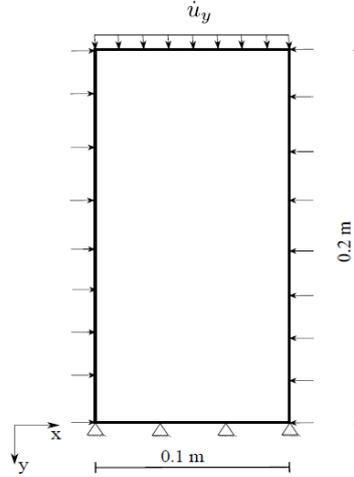

Figure 24: Model setup for example 3

We first present the results of the base simulation in Figures 25 - 30. Figure 25 plots the loading curve on the top boundary of the sample. Figure 26 plots the contours of the equivalent plastic shear strain on the deformed configuration at three displacements, i.e., $u_y = 3 \times 10^{-3}$ m, 4.5 × $10^{-3}$ m, and 6 × $10^{-3}$ m. Figure 27 presents the contours of the plastic volumetric strain at the same loading stages. The positive sign of plastic volume strain denotes dilatation. Figure 28 plots the contours of the micro-rotation at the three loading stages. The results in Figures 26, 27, and 28 have demonstrated the formation of two conjugate shear bands in the sample and the shear band instability nucleates from the sample's geometrical center. It is noted that no weak element, as in the classical finite element method, is needed to trigger the shear banding due to the strong nonlocal formulation in this study [40]. As shown in Figure 28, the micro rotation of material points is localized within the shear band. Figure 29 shows the snapshots of water pressure on the deformed configuration at the three loading stages.

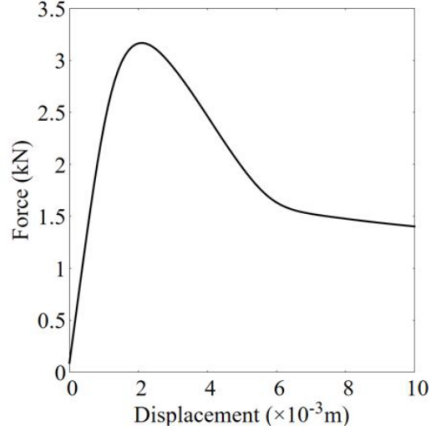

Figure 25: Loading curve on the top boundary.

In this study, the second-order work formulated in the framework of μPPM (see Section 42) is used to validate the shear band zones. Based on the second-order work criterion, the material instability (e.g., shear banding) occurs when the second-order work becomes non-positive. We note that the non-positive second-order work is necessary but not sufficient conditions for shear banding formation in that other types of instability could happen under this condition (e.g., diffusive instability). Figure 30 shows the contours of the second-order work on the deformed configuration. Figure 30 shows that the second-order work within the shear band is negative. It is implied from this that the nonlocal second-order work criterion is useful for detecting the formation of shear bands. We note that in this study, the nonlocal length scale - the horizon is assumed to be equal to the micro-polar length scale. Therefore, the nonlocal length scale in μPPM has a meaningful physical length scale. This could explain the consistency of the results between Figure 26 and Figure 28. In what follows, we study the influence of the spatial discretization schemes, the dilation angle, and the micro-polar length scale on the shear banding instability in this example.

*4.3.1. Influence of spatial discretization*

In this part, we investigate the influence of spatial discretization on the results under the same conditions. For this purpose, we consider two spatial discretization schemes, i.e., 60 × 120 points with $\Delta x = 1.67 \times 10^{-3}$ m and 80 × 160 points with $\Delta x = 1.25 \times 10^{-3}$ m. The same horizon $\delta = 5 \times 10^{-3}$ m is chosen for both cases. All other conditions assumed for the base simulation remain the same. Figure 31 compares the loading curves of the simulations with the two spatial discretization schemes. The results show that the two loading curves are almost identical due to the strong nonlocal formulation. Figure 32 compares the contours of the equivalent plastic shear strain at the displacement $u_y = 6 \times 10^{-3}$ m for the two simulations. Figure 32 compares the contours of plastic volumetric strain at the same displacement. Figure 34 compares the contours of micro rotations at the displacement $u_y = 6 \times 10^{-3}$ m for the two simulations. Figure 35 compares the contour of water pressure at $u_y = 6 \times 10^{-3}$ m. Figure 36 compares the contours of the second- order work at displacement $u_y = 6 \times 10^{-3}$ m for two simulations. The results in Figures 32 - 36 have demonstrated that the numerical results are not dependent on the spatial discretization scheme assuming the same nonlocal length scale due to the strong nonlocal formulations. Next, we present the results of the influence of the micro-polar length on shear banding.

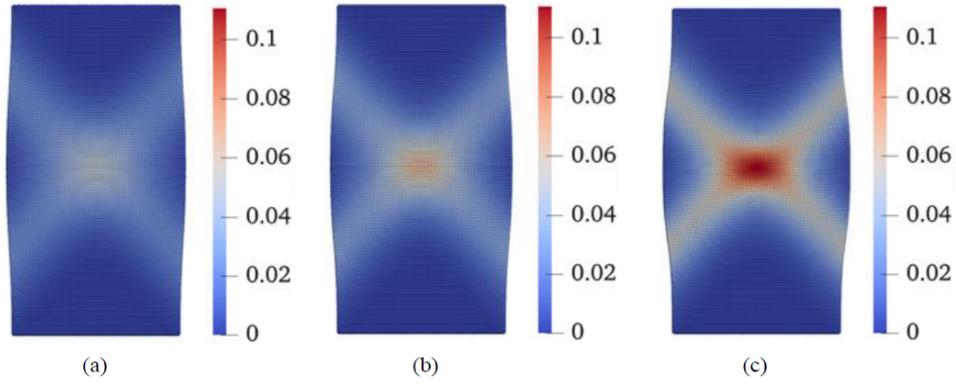

Figure 26: Contours of the equivalent plastic shear strain on the deformed configuration at (a) $u_{y,1} = 3 \times 10^{-3}$ m, (b) $u_{y,2} = 4.5 \times 10^{-3}$ m, and (c) $u_{y,3} = 6 \times 10^{-3}$ m.

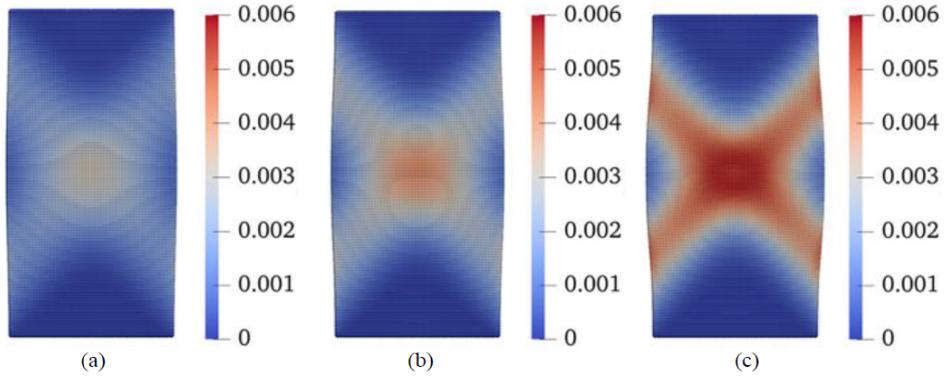

Figure 27: Contours of the plastic volumetric strain on the deformed configuration at (a) $u_{y,1} = 3 \times 10^{-3}$ m, (b) $u_{y,2} = 4.5 \times 10^{-3}$ m, and (c) $u_{y,3} = 6 \times 10^{-3}$ m.

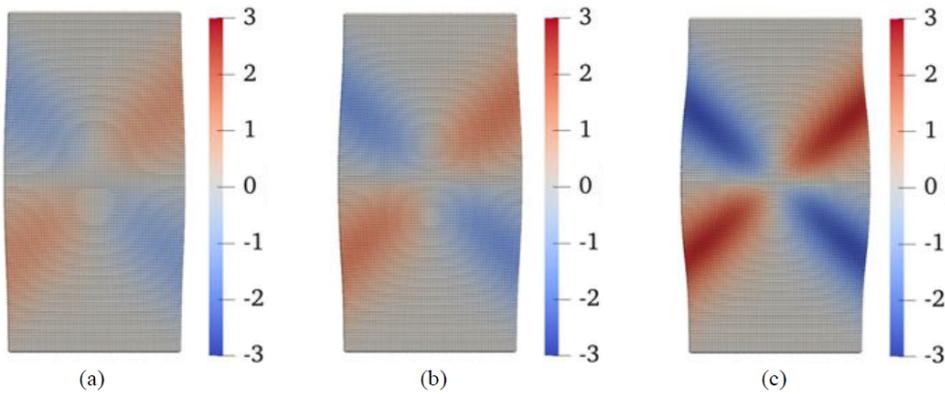

Figure 28: Contours of the micro rotation (degree) on the deformed configuration at (a) $u_{y,1} = 3 \times 10^{-3}$ m, (b) $u_{y,2} = 4.5 \times 10^{-3}$ m, and (c) $u_{y,3} = 6 \times 10^{-3}$ m.

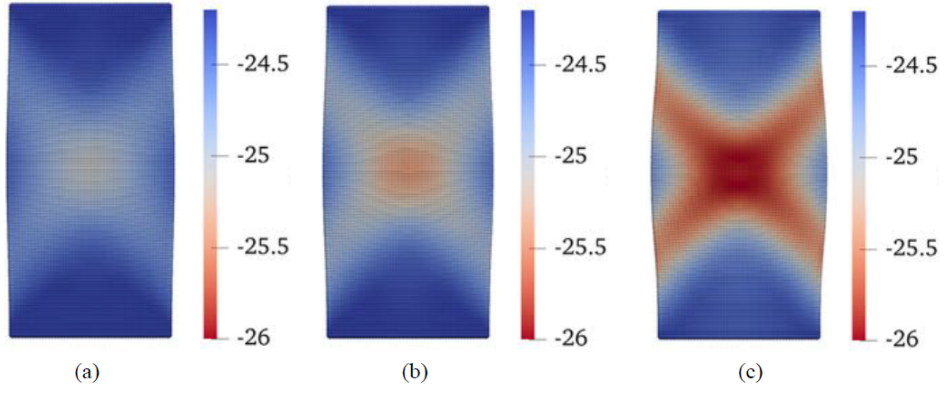

Figure 29: Contours of water pressure (kPa) on the deformed configuration at (a) $u_{y,1} = 3 \times 10^{-3}$ m, (b) $u_{y,2} = 4.5 \times 10^{-3}$ m, and (c) $u_{y,3} = 6 \times 10^{-3}$ m.

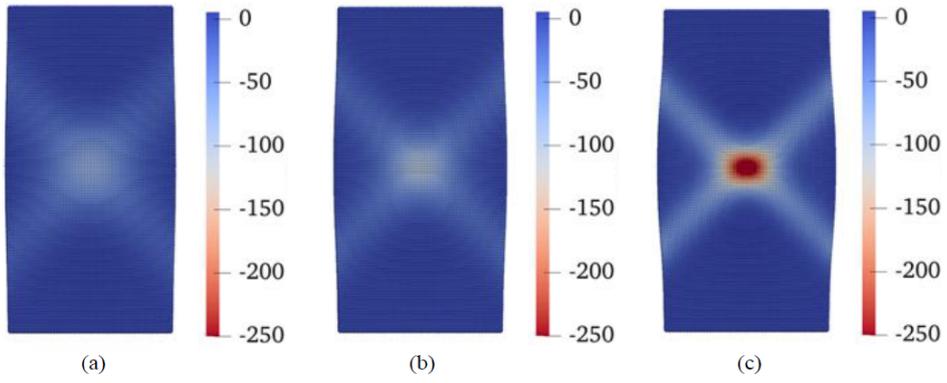

Figure 30: Contours of the second order work (N/m$^2$) on the deformed configuration at (a) $u_{y,1} = 3 \times 10^{-3}$ m, (b) $u_{y,2} = 4.5 \times 10^{-3}$ m, and (c) $u_{y,3} = 6 \times 10^{-3}$ m.

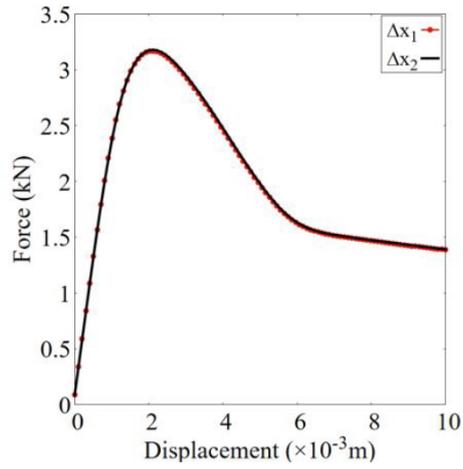

Figure 31: Comparison of the loading curves from the simulations with two spatial discretizations.

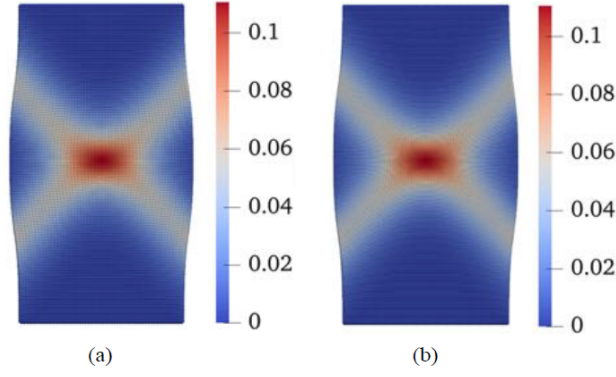

Figure 32: Contours of the equivalent plastic shear strain on the deformed configuration at $u_y = 6 \times 10^{-3}$ m: (a) $\Delta x_1 = 1.67 \times 10^{-3}$ m, and (b) $\Delta x_2 = 1.25 \times 10^{-3}$ m.

*4.3.2. Influence of the micro-polar length scale*

In this part, we evaluate the influence of the micropolar length scale on the shear band formation. For this purpose, we adopt three micropolar length scales, i.e., $l = 3.3 \times 10^{-3}$ m, $5 \times 10^{-3}$ m, and $6.6 \times 10^{-3}$ m. The horizon is the same as the micropolar length scale for each case. The dilatation and frictional angles are the same, i.e., $35°$. All other conditions for the base simulation remain the same. The results are shown in Figures 37 - 42. Figure 37 plots the loading curves on the top boundary for the three simulations. The results show that the micro-polar length scale affects the peak load and the post-localization stage. Specifically, the simulations with a large micro-polar length scale generate a relatively large peak load and a higher post-localization loading curve. We note that in the present study, the horizon equals the micro-polar length scale, which is different from the previous study assuming a horizon without a clear physical meaning [28]. Figure 38 plots the contours of the equivalent plastic shear strain on the deformed configuration at $u_y = 6 \times 10^{-3}$ m from the simulations with three different micropolar length scales, i.e., $l = 3.3 \times 10^{-3}$ m, $5 \times 10^{-3}$ m, and $6.6 \times 10^{-3}$ m. Figure 39 presents contours of the plastic volumetric strain in the three simulations at $u_y = 6 \times 10^{-3}$ m. Figure 40 plots the contours of micro rotation for the three simulations. Figure 41 compares the contours of water pressure for the three simulations at the same loading stage.

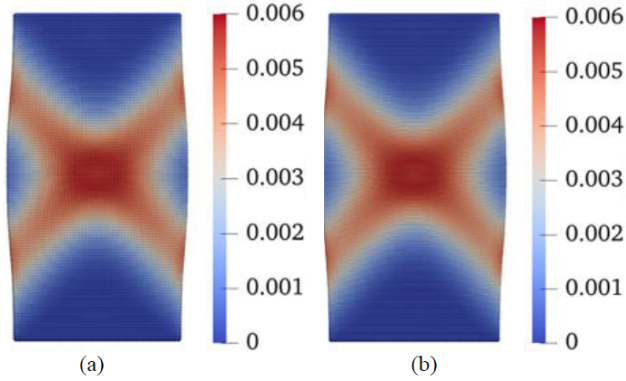

Figure 33: Contours of the plastic volumetric strain on the deformed configuration at $u_y = 6 \times 10^{-3}$ m: (a) $\Delta x_1 = 1.67 \times 10^{-3}$ m, and (b) $\Delta x_2 = 1.25 \times 10^{-3}$ m.

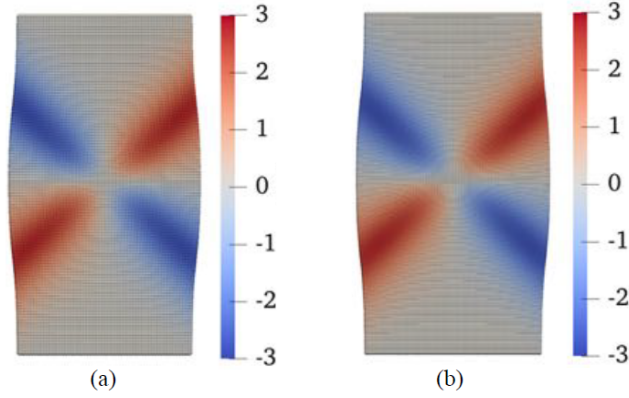

Figure 34: Contours of micro rotation (degree) on the deformed configuration at $u_y = 6 \times 10^{-3}$ m for the cases: (a) $\Delta x_1 = 1.67 \times 10^{-3}$ m, and (b) $\Delta x_2 = 1.25 \times 10^{-3}$ m.

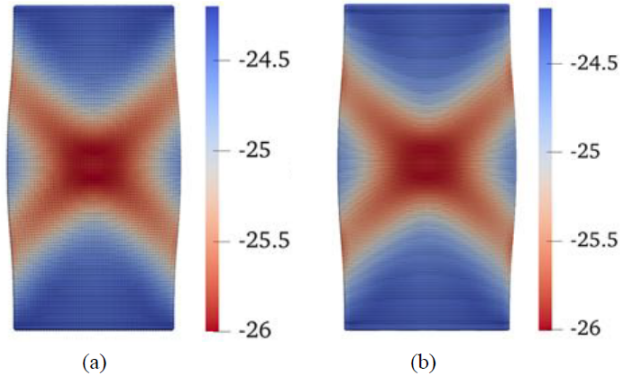

Figure 35: Contours of water pressure (kPa) on the deformed configuration at $u_y = 6 \times 10^{-3}$ m: (a) $\Delta x_1 = 1.67 \times 10^{-3}$ m, and (b) $\Delta x_2 = 1.25 \times 10^{-3}$ m.

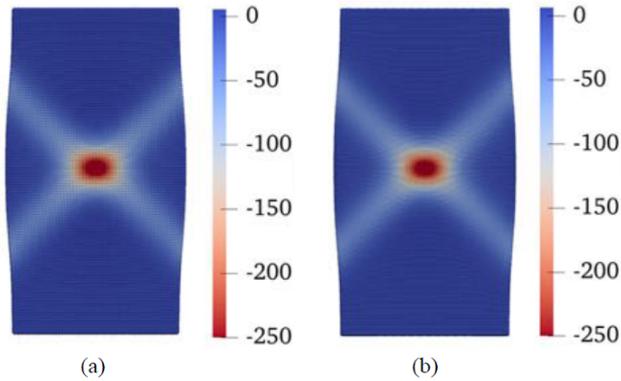

Figure 36: Contours of the second order work (N/m$^2$) on the deformed configuration at $u_y = 6 \times 10^{-3}$ m: (a) $\Delta x_1 = 1.67 \times 10^{-3}$ m, and (b) $\Delta x_2 = 1.25 \times 10^{-3}$ m.

The results in Figures 38, 39, 40, and 41 demonstrate that the micro-polar length scale affects the shear band width but not the shear band inclination. Specifically, the shear band width may decrease in the numerical simulations by decreasing the micro-polar length scale. It is implied from this finding that the experimental data on shear band width can be used to calibrate the micro-length scale and the nonlocality of the µPPM. Comparison of the results in Figures 38 and 40 demonstrates that the shear band width manifested in both the equivalent plastic shear strain and the micro-rotation of material points are consistent. What is implied from this finding is that the micro-rotation of material points is concentrated into the banded zoned, as found in the laboratory testing of shear bands [41]. Figure 42 shows that the second-order work on the deformed configuration for the three simulations at $u_y = 6 \times 10^{-3}$ m. The negative second-order work within the shear band for all three cases corroborates the shear band instability in the sample.

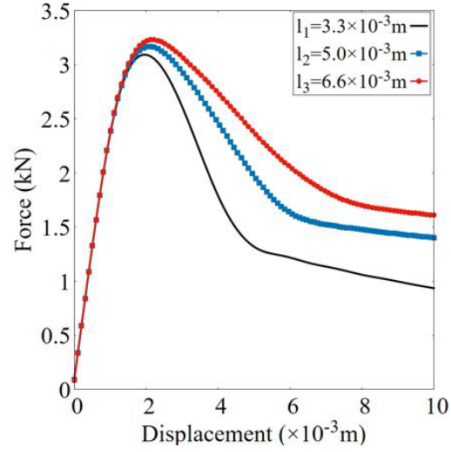

Figure 37: Loading curves on the top boundary from the simulations with the three micropolar length scales, $l = 3.3 \times 10^{-3}$ m, $5 \times 10^{-3}$ m, and $6.6 \times 10^{-3}$ m

### 4.3.3. Influence of the dilation angle

It is known that the dilation angle could impact the shear band angle. Thus, in this part, we investigate the influence of the dilatation angle on the inclination angle of the shear band through the coupled µPPM. For this purpose, we run simulations with three dilatation angles, i.e., 15°, 25°, and 35°. The same frictional angle 35° is assumed for the three simulations. All the other conditions for the base simulation are the same. The results are reported in Figures 43 - 48.

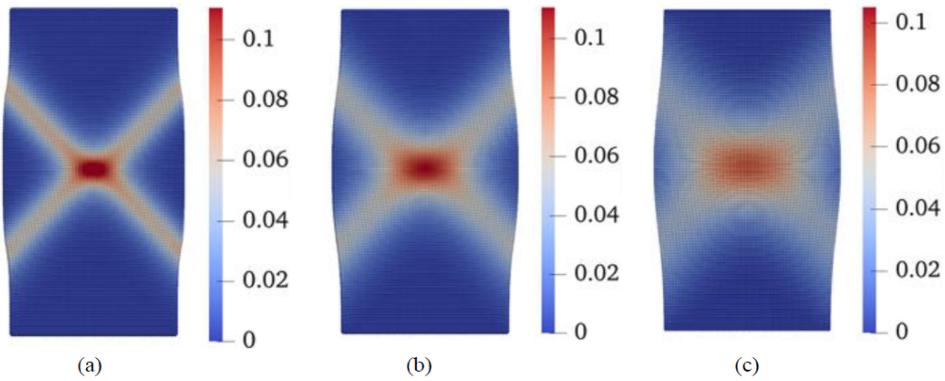

Figure 38: Contours of the equivalent shear plastic strain on the deformed configuration at $u_y = 6 \times 10^{-3}$ m from the simulations: (a) $l_1 = 3.3 \times 10^{-3}$ m, (b) $l_2 = 5 \times 10^{-3}$ m, and (c) $l_3 = 6.6 \times 10^{-3}$ m.

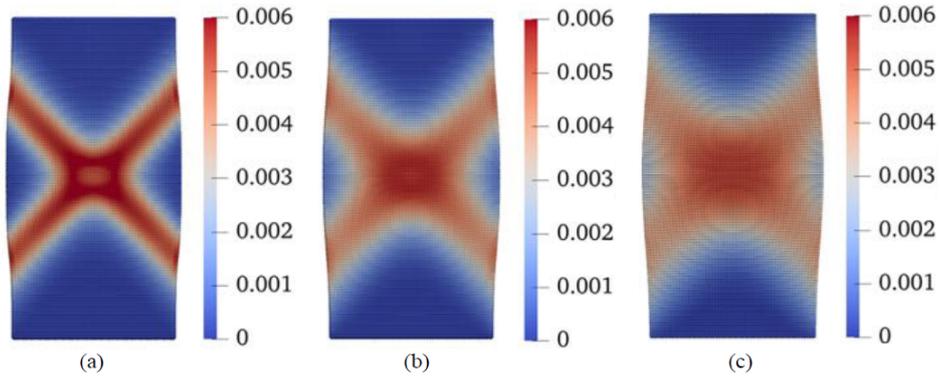

Figure 39: Contours of the plastic volume strain on the deformed configuration at $u_y = 6 \times 10^{-3}$ m from the simulations: (a) $l_1 = 3.3 \times 10^{-3}$ m, (b) $l_2 = 5 \times 10^{-3}$ m, and (c) $l_3 = 6.6 \times 10^{-3}$ m.

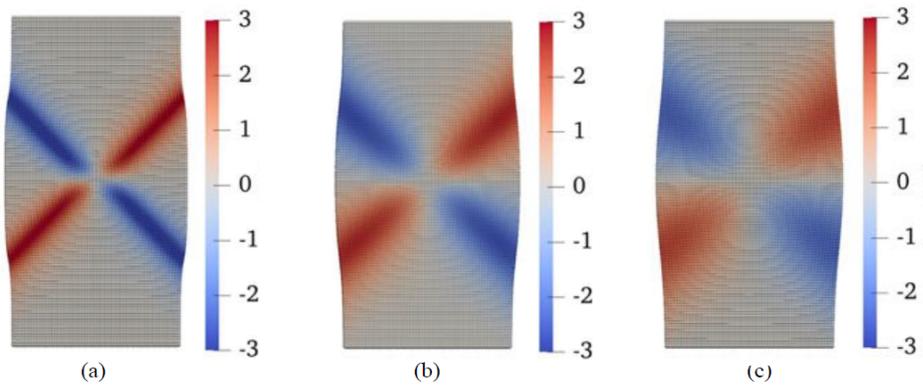

Figure 40: Contours of the micro rotation (degree) on the deformed configuration at $u_y = 6 \times 10^{-3}$ m from the simulations: (a) $l_1 = 3.3 \times 10^{-3}$ m, (b) $l_2 = 5 \times 10^{-3}$ m, and (c) $l_3 = 6.6 \times 10^{-3}$ m.

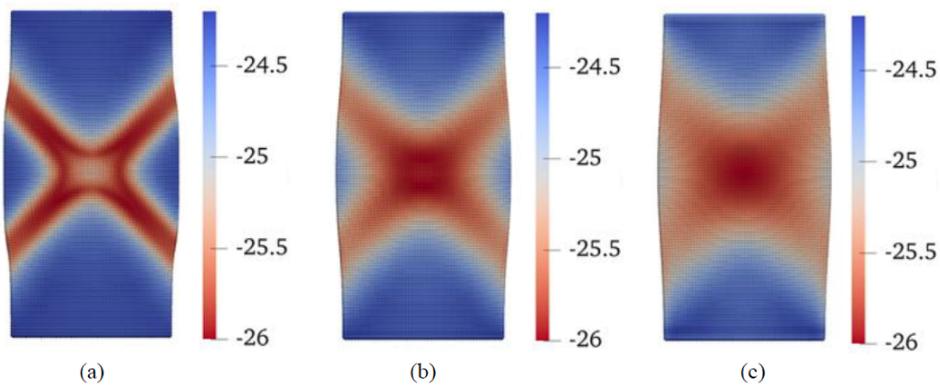

Figure 41: Contours of water pressure (kPa) on the deformed configuration at $u_y = 6 \times 10^{-3}$ m from the simulations: (a) $l_1 = 3.3 \times 10^{-3}$ m, (b) $l_2 = 5 \times 10^{-3}$ m, and (c) $l_3 = 6.6 \times 10^{-3}$ m.

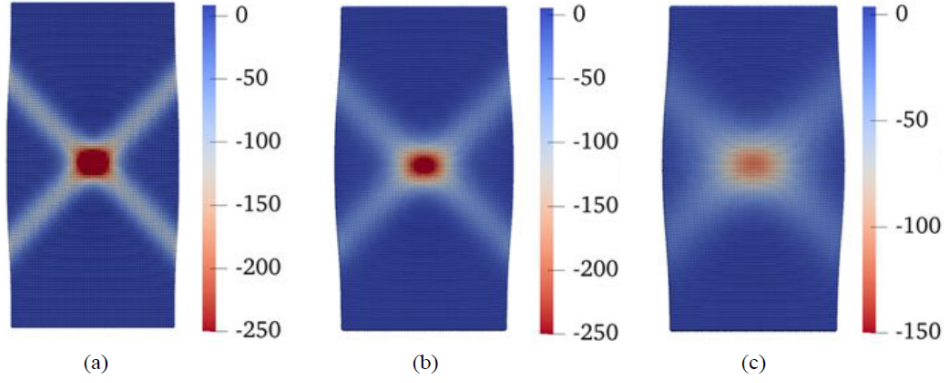

Figure 42: Contours of the second order work (N/m$^2$) on the deformed configuration at $u_y = 6 \times 10^{-3}$ m from the simulations: (a) $l_1 = 3.3 \times 10^{-3}$ m, (b) $l_2 = 5 \times 10^{-3}$ m, and (c) $l_3 = 6.6 \times 10^{-3}$ m.

Figure 43 plots the loading curve on the top boundary for the three simulations with different dilation angles. The results in Figure 43 imply that the dilation angle does not affect the pre-localization load curve and the peak load, while it may slightly affect the post-localization load curve. Figure 44 plots the contours of the equivalent plastic shear strain on the deformed configuration at at $u_y = 6 \times 10^{-3}$ m from the simulations with three different dilation angles. Similarly, Figure 45 presents the contours of the plastic volumetric strain on the deformed configuration at displacement at $u_y = 6 \times 10^{-3}$ m for the three simulations. Figure 46 plots the contours of micro rotation on the deformed configuration at displacement at $u_y = 6 \times 10^{-3}$ m. Figure 47 shows the water pressure contours on the deformed configuration at $u_y = 6 \times 10^{-3}$ m for the three simulations. Figure 48 plots the second-order work contour on the deformed configuration at displacement at $u_y = 6 \times 10^{-3}$ m. Next, we briefly discuss the impact of dilation angles on the shear banding formation, as implied by the results in the aforementioned figures. First, the results in Figures 44-48 show that the shear band inclination angle can be affected by the dilation angle. Specifically, the shear band inclination angle with respect to the horizontal direction increases with the dilation angle increase. This observation agrees with Roscoe's solution for the shear banding at a material point, i.e., a homogeneous specimen (e.g., [41, 51]) Second, as shown in Figure 45, the dilatation angle affects the plastic volumetric strain and variation of fluid pressure in the development of shear bands. For example, for the simulation with the dilatation angle 35°, the plastic volumetric strain is positive (i.e., dilative), while the plastic volumetric strain from the simulations with 15° and 25° are negative (i.e., compaction). Accordingly, the results in Figure 47 show that the matric suction (i.e., negative fluid pressure) in the banded zone increases for the simulation with 35° due to the solid skeleton dilation while the matric suction in shear bands from the simulations with 15° and 25° decreases due to the solid skeleton compaction. In summary, these findings show that the dilation angle of the porous media affects the shear band angle and impacts the fluid pressure in shear bands under certain loading conditions. In the following section, we present the numerical simulation of a forward-running crack in unsaturated porous media.

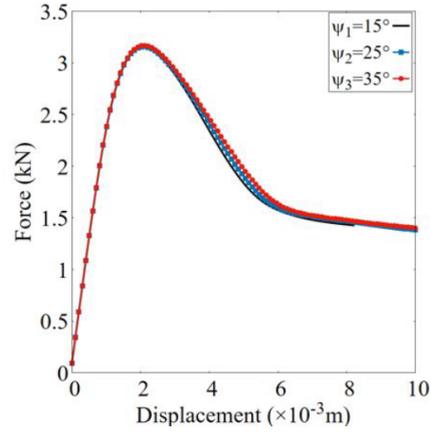

Figure 43: Loading curves on the top boundary from the simulations with $\psi = 15°$, $25°$ and $35°$.

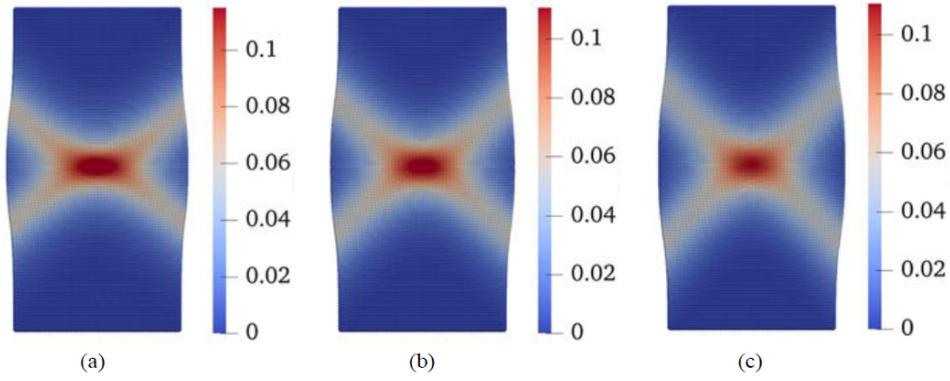

Figure 44: Contours of the equivalent shear plastic strain on the deformed configuration at $u_y = 6 \times 10^{-3}$ m from the simulations with dilation angle (a) $\psi_1 = 15°$, (b) $\psi_2 = 25°$, and (c) $\psi_3 = 35°$.

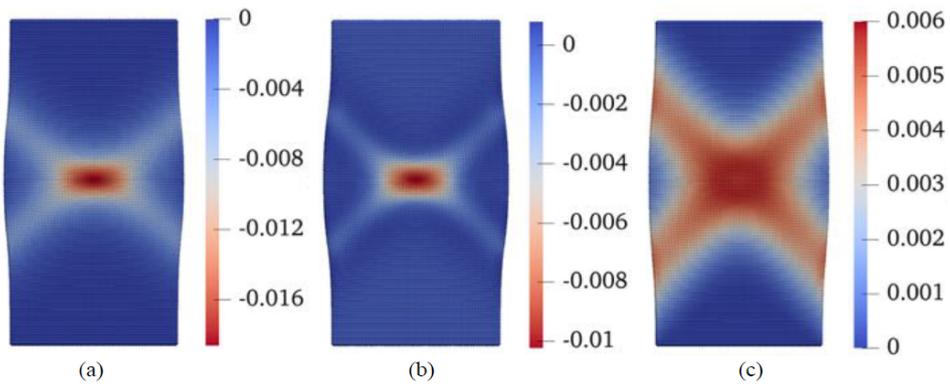

Figure 45: Contours of plastic volumetric strain on the deformed configuration at $u_y = 6 \times 10^{-3}$ m: (a) $\psi_1 = 15°$, (b) $\psi_2 = 25°$, and (c) $\psi_3 = 35°$.

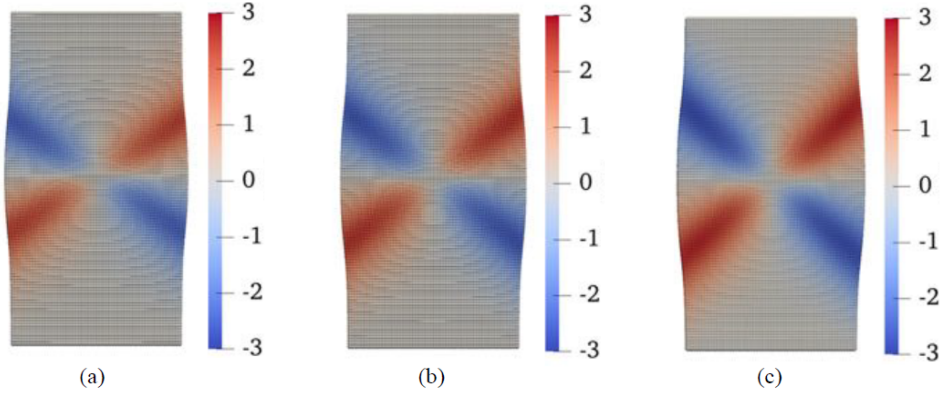

Figure 46: Contours of micro rotation (degree) on the deformed configuration at $u_y = 6 \times 10^{-3}$ m: (a) $\psi_1 = 15°$, (b) $\psi_2 = 25°$, and (c) $\psi_3 = 35°$.

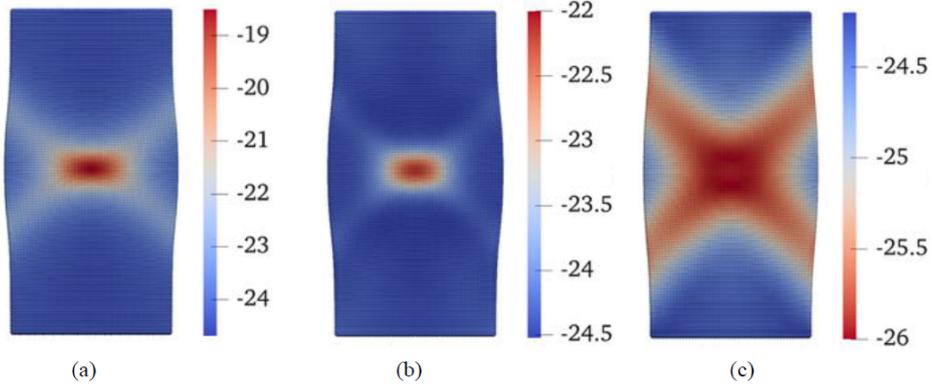

Figure 47: Contours of water pressure (kPa) on the deformed configuration at $u_y = 6 \times 10^{-3}$ m: (a) $\psi_1 = 15°$, (b) $\psi_2 = 25°$, and (c) $\psi_3 = 35°$.

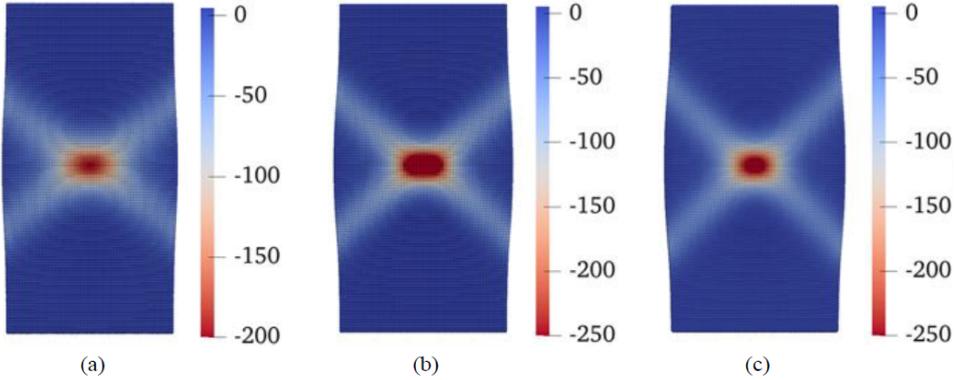

Figure 48: Contours of the second order work (N/m$^2$) on the deformed configuration at $u_y = 6 \times 10^{-3}$ m: (a) $\psi_1 = 15°$, (b) $\psi_2 = 25°$, and (c) $\psi_3 = 35°$.

### 4.4. Example 4: Forward running crack propagation

This example deals with forward-running crack propagation in unsaturated elastic porous material. The proposed micro-polar J-integral bond-breakage criterion is adopted to model crack propagation. In this example, we will study the impact of micro-polar length scales on crack growth. Figure 49 plots the geometry, the initial crack location, and the loading condition for this example. The initial horizontal crack along the center line of the sample is 0.05 m in length, as shown in Figure 49. The vertical displacement load of $u_y = 6 \times 10^{-5}$ m is imposed on the top and bottom boundaries. The left and right boundaries are free to deform. All fluid boundaries are impervious. The initial water pressure is $p_0 = 1$ MPa. The material parameters are: solid density $\rho_s$ = 2650 kg/m3, initial porosity $\phi_0$ = 0.35, bulk modulus K = 16.67 GPa, shear modulus $\mu$ = 10 GPa, micropolar shear modulus $\mu c$ = 5 GPa, micropolar length scale = $l = 2 \times 10^{-3}$ m, water viscosity $\mu w$ = $1 \times 10^{-9}$ MPa, water density $\rho_w$ = 1000 kg/m, hydraulic conductivity kw= $1 \times 10^{-8}$ m/s, n = 1.8 and sa = 10 MPa. For this example, Gcr = 100 N/m is adopted for the J integral bond breakage criterion. The stabilization parameter G = 0.5 is assumed. The domain is discretized to 200 × 150 material points with 6.x = 5 ⇁ 10 m. The horizon $\delta$ = 4.05Δ.x. The simulation time $t = 6.25 \times 10^{-5}$ s, and the time increment is $\Delta t = 2.5 \times 10^{-8}$ s.

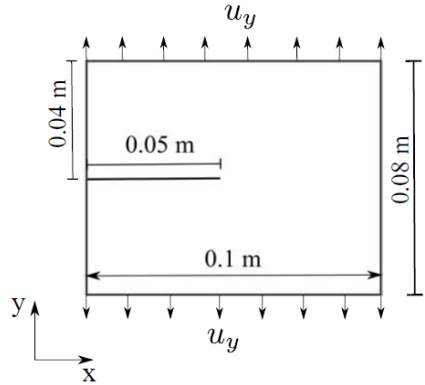

Figure 49: Model setup for example 4.

First, we present the results of a base simulation. Figure 50 plots the snapshots of the damage variable on the deformed configuration at three loading steps, i.e., $u_y$ = $0.9 \times 10^{-5}$ m, $1.2 \times 10^{-5}$ m, and $1.5 \times 10^{-5}$ m. As shown in Figure 50, the crack growth starts at $u_y$ = $0.9 \times 10^{-5}$ m and the crack length reaches 0.078 m at uy = $u_y$ = $1.5 \times 10^{-5}$ m. Figure 51 presents the contour of water pressure at these three loading stages. It is implied from Figure 51 that the crack propagation increases the water pressure at the crack tip and decreases the water pressure at the crack path. This result is due to the water flow from bulk to fracture space during crack propagation. Figure 52 presents the contours of the micro rotation of material points at these three loading stages. Figure 52 shows that the micro rotation of material points is concentrated on the crack tip as the crack propagates forward. The magnitude of micro rotations increases with the crack growth.

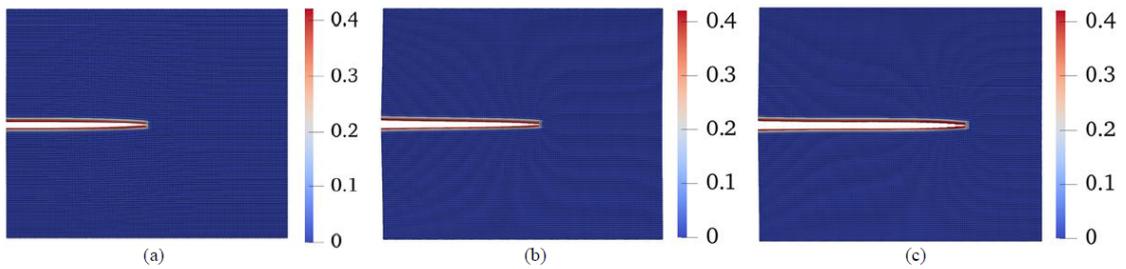

Figure 50: Contours of the damage variable on the deformed configuration (MF = 100) at three loading stages: (a) $u_{y,1} = 0.9 \times 10^{-5}$ m, (b) $u_{y,2} = 1.2 \times 10^{-5}$ m, and (c) $u_{y,3} = 1.5 \times 10^{-5}$ m.

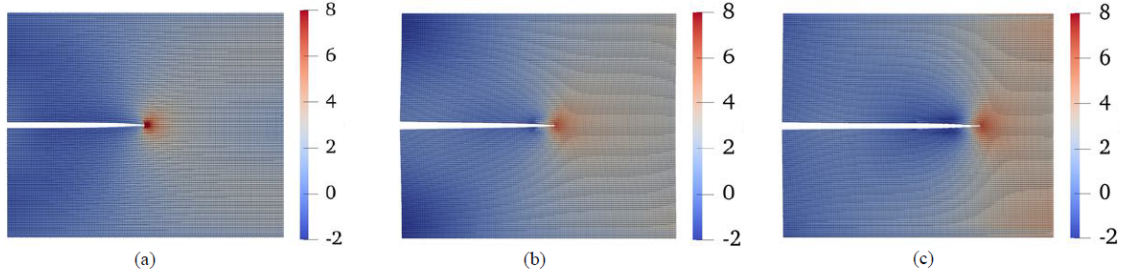

Figure 51: Contours of water pressure (MPa) on the deformed configuration (MF = 100) at three loading stages: (a) $u_{y,1} = 0.9 \times 10^{-5}$ m, (b) $u_{y,2} = 1.2 \times 10^{-5}$ m, and (c) $u_{y,3} = 1.5 \times 10^{-5}$ m.

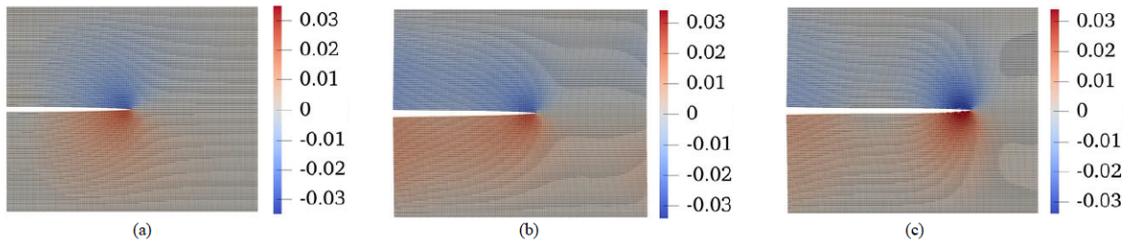

Figure 52: Contours of micro rotation (degree) on the deformed configuration (MF = 100) at three loading stages: (a) $u_{y,1} = 0.9 \times 10^{-5}$ m, (b) $u_{y,2} = 1.2 \times 10^{-5}$ m, and (c) $u_{y,3} = 1.5 \times 10^{-5}$ m.

Second, we study the sensitivity of the numerical results to the spatial discretization schemes in this example. For this purpose, we present the results of the simulations with two uniform grids. The same micro-polar and horizon are adopted, i.e., $2 \times 10^{-3}$ m. For the two cases, the specimen is discretized into $150 \times 120$ uniform material points ($\Delta x_1 = 6.7 \times 10^{-3}$ m) and $200 \times 150$ uniform material points ($\Delta x_2 = 6.7 \times 10^{-3}$ m), respectively. The other material parameters remain the same. The results are presented in Figures 53, 54, and 55. Figure 53 plots the contours of the damage variable on the deformed configuration at the same loading stage, i.e., $u_y = 1.5 \times 10^{-5}$ m, for both simulations. Figure 54 shows the contours of water pressure at $u_y = 1.5 \times 10^{-5}$ m for the two simulations. Figure 55 plots the contours of the micro rotation at $u_y = 1.5 \times 10^{-5}$ m for the two simulations. The results have demonstrated that with the same micro-polar length scale, the numerical results are insensitive to the spatial discretization scheme.

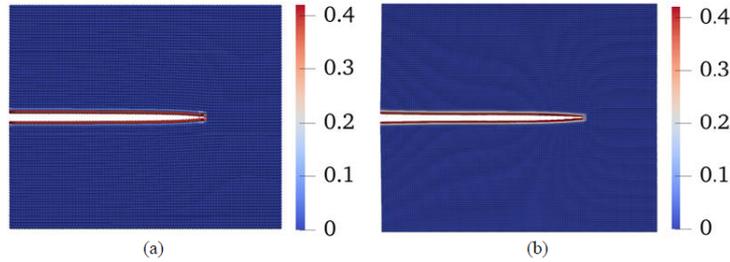

Figure 53: Contours of the damage variable on the deformed configuration (MF = 100) at the same loading stage (i.e., $u_y = 1.5 \times 10^{-5}$ m): (a) $\Delta x_1 = 6.7 \times 10^{-4}$ m and (b) $\Delta x_2 = 5 \times 10^{-4}$ m.

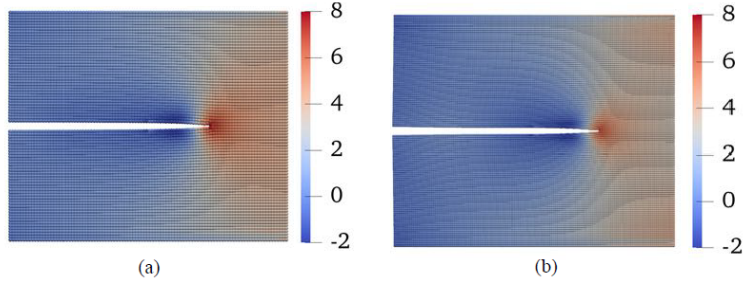

Figure 54: Contours of water pressure (MPa) on the deformed configuration (MF $= 100$) at $u_y = 1.5 \times 10^{-5}$ m: (a) $\Delta x_1 = 6.7 \times 10^{-4}$ m and (b) $\Delta x_2 = 5 \times 10^{-4}$ m.

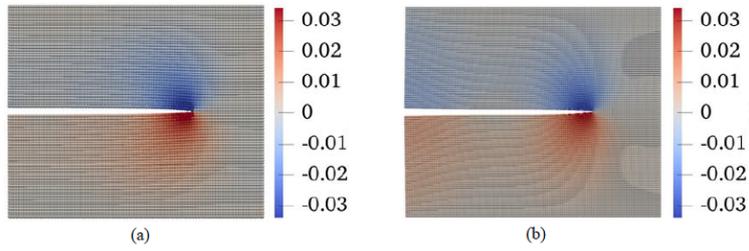

Figure 55: Contours of micro rotation (degree) on the deformed configuration (MF $= 100$) at $u_y = 1.5 \times 10^{-5}$ m: (a) $\Delta x_1 = 6.7 \times 10^{-4}$ m and (b) $\Delta x_2 = 5 \times 10^{-4}$ m.

Third, we investigate the impact of the micropolar length scale on the J integral. For this purpose, we compare the numerical results assuming three micro-polar length scales, i.e., $l = 1 \times 10^{-3}$ m, $2 \times 10^{-3}$ m and $4 \times 10^{-3}$ m. The other parameters and conditions are the same as the base simulation. The horizon equals the micropolar length scale. Figure 56 plots the contours of the damage variable on the deformed configuration at $u_y = 1.5 \times 10^{-5}$ m for the simulations with three micro-polar length scales. As shown in Figure 56 at the same displacement load $u_y = 1.5 \times 10^{-5}$ m the crack length reaches 0.08 m, 0.078 m, and 0.075 m for the three simulations with $l = 1 \times 10^{-3}$, $l = 2 \times 10^{-3}$, and $l = 4 \times 10^{-3}$, respectively. Figure shows the contours of water pressure at $u_y = 1.5 \times 10^{-5}$ m for the three simulations. Figure plots the contours of the micro rotation on the deformed configuration at $u_y = 1.5 \times 10^{-5}$ m. Figures 57 and 58 show that the water pressure and micro-rotation at the crack tip increase by decreasing the micropolar length scale. Figure 59 plots the variation of the J integral with a/w from the three simulations. Here a is the crack length, and w is the width of the specimen. Figure 59 shows that the micro-polar length scale affects the J integral in this example. The simulation with a larger micro-polar length scale has a larger value of the J-integral at the same loading stage. This finding may imply that the micro-polar length scale can be determined from the crack energy release rate as an inverse problem given the latter from the physical testing.

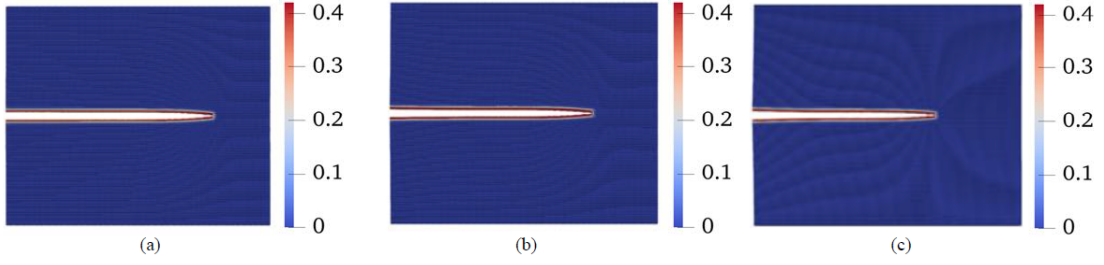

Figure 56: Contours of the damage variable on the deformed configuration (MF = 100) at $u_y = 1.5 \times 10^{-5}$ m: (a) $l_1 = 1 \times 10^{-3}$ m, (b) $l_2 = 2 \times 10^{-3}$ m, and (c) $l_3 = 4 \times 10^{-3}$ m.

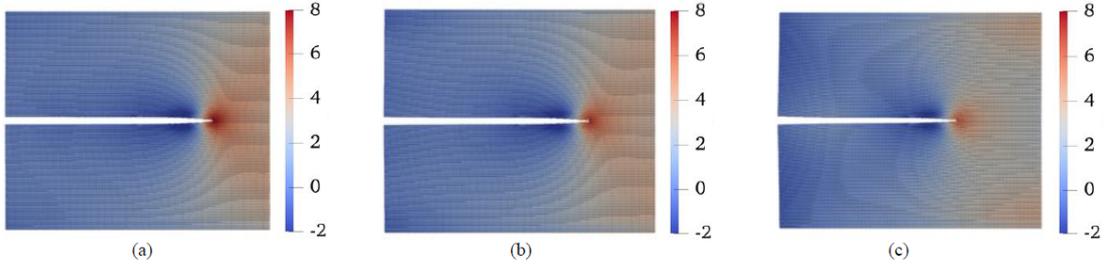

Figure 57: Contours of water pressure (MPa) on the deformed configuration (MF = 100) at $u_y = 1.5 \times 10^{-5}$ m: (a) $l_1 = 1 \times 10^{-3}$ m, (b) $l_2 = 2 \times 10^{-3}$ m, and (c) $l_3 = 4 \times 10^{-3}$ m.

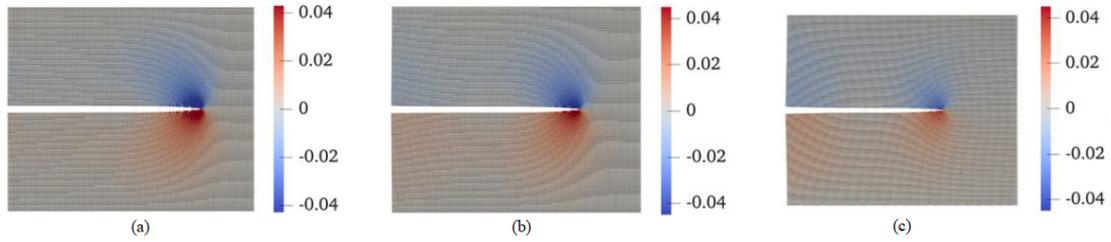

Figure 58: Contours of micro rotation (degree) on the deformed configuration (MF = 100) at $u_y = 1.5 \times 10^{-5}$ m: (a) $l_1 = 1 \times 10^{-3}$ m, (b) $l_2 = 2 \times 10^{-3}$ m, and (c) $l_3 = 4 \times 10^{-3}$ m.

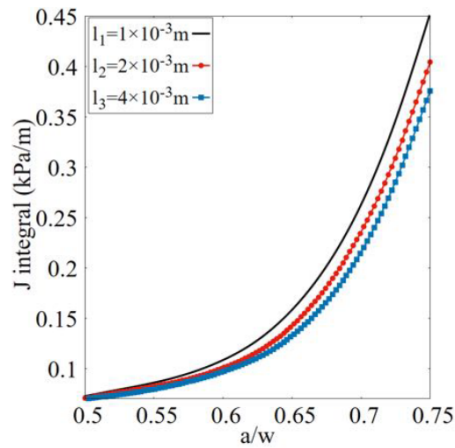

Figure 59: Comparison of the J integral-$a/w$ curves for the three simulations with $l = 1 \times 10^{-3}$ m, $2 \times 10^{-3}$ m, and $4 \times 10^{-3}$ m.

## 5. Closure

In this article, we have developed a multiphase μPPM paradigm for modeling dynamic shear banding and fracturing in unsaturated porous media. As a novelty, in the coupled μPPM, a microstructure-based material length scale is incorporated by considering micro-rotations of the solid skeleton in line with the Cosserat continuum theory for solids. As a new contribution, we reformulate the second-order work for detecting material instability and the energy-based crack criterion and J-integral for modeling fracturing in the μPPM paradigm. The stabilized Cosserat PPM correspondence principle augmented to include unsaturated fluid flow is adopted to circum- vent the multiphase zero-energy mode instability in the proposed μPPM. We have numerically implemented the novel μPPM paradigm through a fractional-step explicit-explicit algorithm in time that splits the coupled problem into a deformation/fracture problem and an unsaturated fluid flow problem. A hybrid Lagrangian-Eulerian meshfree method is utilized to discretize the problem domain in space. Numerical examples are presented to evaluate the numerical implementation and demonstrate the efficacy of the proposed μPPM paradigm for modeling shear banding and fracturing in unsaturated porous media. The numerical results have shown that this fully coupled μPPM paradigm is robust to model shear band instability and fracturing in variably saturated porous media.

## Acknowledgements

This work has been supported by the US National Science Foundation under contract numbers 1659932 and 1944009. The support is gratefully acknowledged. Any opinions or positions expressed in this article are those of the authors only and do not reflect any opinions or positions of the NSF.

## References

[1] Zienkiewicz OC, Chan A, Pastor M, Schrefler B, Shiomi T. Computational geomechanics; vol. 613. Citeseer; 1999.
[2] Lewis RW, Schrefler B. The finite element method in the static and dynamic deformation and consolidation of porous media. John Wiley & Sons; 1998.
[3] Fung Yc. Biomechanics: mechanical properties of living tissues. Springer Science & Business Media; 2013.
[4] Cowin SC, Doty SB. Tissue mechanics. Springer; 2007.
[5] Terzaghi K, Peck RB. Soil mechanics. Engineering Practice John Wiley and Sons, Inc, New York 1948;.
[6] Fredlund DG, Rahardjo H, Fredlund MD. Unsaturated Soil Mechanics in Engineering Practice. John Wiley & Sons; 2012.
[7] de Borst R. Computational methods for fracture in porous media: Isogeometric and extended finite element methods. Elsevier; 2017.
[8] Wang K, Song X. Strain localization in non-isothermal unsaturated porous media considering material heterogeneity with stabilized mixed finite elements. Computer Methods in Applied Mechanics and Engineering 2020;359:112770.
[9] Song X, Wang K, Ye M. Localized failure in unsaturated soils under non-isothermal conditions. Acta Geotechnica 2018;13:73–85.
[10] Song X. Transient bifurcation condition of partially saturated porous media at finite strain. International Journal for Numerical and Analytical Methods in Geomechanics 2017;41(1):135–56.
[11] Song X, Ye M, Wang K. Strain localization in a solid-water-air system with random heterogeneity via stabilized mixed finite elements. International Journal for Numerical Methods in Engineering 2017;112(13):1926–50.
[12] Song X, Borja RI. Mathematical framework for unsaturated flow in the finite deformation range. International Journal for Numerical Methods in Engineering 2014;97(9):658–82.
[13] Seed HB, Lee KL, Idriss IM, Makdisi FI. The slides in the san fernando dams during the earthquake of february 9, 1971. Journal of the Geotechnical Engineering Division 1975;101(7):651–88.
[14] Bray JD, Seed RB, Cluꜛ LS, Seed HB. Earthquake fault rupture propagation through soil. Journal of Geotechnical Engineering 1994;120(3):543–61.
[15] Kramer SL. Geotechnical earthquake engineering. Pearson Education India; 1996.
[16] Bray JD, Rodriguez-Marek A. Characterization of forward-directivity ground motions in the near-fault region. Soil dynamics and earthquake engineering 2004;24(11):815–28.
[17] Matasovic N, Kavazanjian Jr E, Augello AJ, Bray JD, Seed RB. Solid waste landfill damage caused by 17 january 1994 northridge earthquake. The Northridge, California, Earthquake of 1994;17:221–9.
[18] Bray JD, Travasarou T. Simplified procedure for estimating earthquake-induced deviatoric slope displacements. Journal of geotechnical and geoenvironmental engineering 2007;133(4):381–92.
[19] Wartman J, Seed RB, Bray JD. Shaking table modeling of seismically induced deformations in slopes. Journal of Geotechnical and Geoenvironmental Engineering 2005;131(5):610–22.
[20] Wartman J, Montgomery DR, Anderson SA, Keaton JR, Benoˆıt J, dela Chapelle J, et al. The 22 march 2014 oso landslide, washington, usa. Geomorphology 2016;253:275–88.
[21] Alonso EE. Triggering and motion of landslides. G´eotechnique 2021;71(1):3–59.

[22] Khoei AR. Extended finite element method: theory and applications. John Wiley & Sons; 2014.
[23] Menon S, Song X. Coupled analysis of desiccation cracking in unsaturated soils through a non-local mathematical formulation. Geosciences 2019;9(10):428.
[24] Song X, Menon S. Modeling of chemo-hydromechanical behavior of unsaturated porous media: a nonlocal approach based on integral equations. Acta Geotechnica 2019;14(3):727–47.
[25] Song X, Khalili N. A peridynamics model for strain localization analysis of geomaterials. International Journal for Numerical and Analytical Methods in Geomechanics 2019;43(1):77–96.
[26] Menon S, Song X. Shear banding in unsaturated geomaterials through a strong nonlocal hydromechanical model. European Journal of Environmental and Civil Engineering 2020;:1–15.
[27] Song X, Silling SA. On the peridynamic effective force state and multiphase constitutive correspondence principle. Journal of the Mechanics and Physics of Solids 2020;145:104161.
[28] Menon S, Song X. A computational periporomechanics model for localized failure in unsaturated porous media. Computer Methods in Applied Mechanics and Engineering 2021;384:113932.
[29] Menon S, Song X. Computational multiphase periporomechanics for unguided cracking in unsaturated porous media. International Journal for Numerical Methods in Engineering 2022;123(12):2837–71.
[30] Menon S, Song X. A stabilized computational nonlocal poromechanics model for dynamic analysis of saturated porous media. International Journal for Numerical Methods in Engineering 2021;122(20):5512–39.
[31] Menon S, Song X. Updated lagrangian unsaturated periporomechanics for extreme large deformation in unsaturated porous media. Computer Methods in Applied Mechanics and Engineering 2022;400:115511.
[32] Menon S, Song X. Computational coupled large-deformation periporomechanics for dynamic failure and fracturing in variably saturated porous media. International Journal for Numerical Methods in Engineering 2022;.
[33] Coussy O. Poromechanics. John Wiley & Sons; 2004.
[34] Cheng AHD. Poroelasticity; vol. 27. Springer; 2016.
[35] Silling SA, Epton M, Weckner O, Xu J, Askari E. Peridynamic states and constitutive modeling. Journal of elasticity 2007;88(2):151–84.
[36] Madenci E, Oterkus E. Peridynamic theory and its applications. Springer; 2014.
[37] Silling SA. Stability of peridynamic correspondence material models and their particle discretizations. Computer Methods in Applied Mechanics and Engineering 2017;322:42–57.
[38] Cosserat EMP, Cosserat F. Théorie des corps déformables. A. Hermann et fils; 1909.
[39] Eringen AC, Eringen AC. Theory of micropolar elasticity. Springer; 1999.
[40] Mühlhaus HB, Vardoulakis I. The thickness of shear bands in granular materials. Geotechnique 1987;37(3):271–83.
[41] Sulem J, Vardoulakis I. Bifurcation analysis in geomechanics. CRC Press; 1995.
[42] de Borst R. A generalisation of j2-flow theory for polar continua. Computer Methods in Applied Mechanics and Engineering 1993;103(3):347–62.
[43] Ehlers W, Volk W. On shear band localization phenomena of liquid-saturated granular elastoplastic porous solid materials accounting for fluid viscosity and micropolar solid rotations. Mechanics of Cohesive-frictional Materials: An International Journal on Experiments, Modelling and Computation of Materials and Structures 1997;2(4):301–20.
[44] Ehlers W, Volk W. On theoretical and numerical methods in the theory of porous media based on polar and non-polar elasto-plastic solid materials. International Journal of Solids and Structures 1998;35(34-35):4597–617.
[45] Mühlhaus H, Pasternak E. Path independent integrals for cosserat continua and application to crack problems. International journal of fracture 2002;113:21–6.
[46] Tang H, Wei W, Song X, Liu F. An anisotropic elastoplastic cosserat continuum model for shear failure in stratified geomaterials. Engineering Geology 2021;293:106304.
[47] Rice JR. A path independent integral and the approximate analysis of strain concentration by notches and cracks 1968;.
[48] Anderson TL. Fracture mechanics: fundamentals and applications. CRC press; 2017.
[49] Gerstle W, Sau N, Silling S. Peridynamic modeling of concrete structures. Nuclear engineering and design 2007;237(12-13):1250–8.
[50] Chowdhury SR, Rahaman MM, Roy D, Sundaram N. A micropolar peridynamic theory in linear elasticity. International Journal of Solids and Structures 2015;59:171–82.
[51] Song X, Pashazad H. Visco-cosserat periporomechanics for dynamic shear bands and crack branching in porous media. Arxiv 2023;:1 – 42.
[52] Silling SA, Lehoucq RB. Peridynamic theory of solid mechanics. Advances in applied mechanics 2010;44:73–168.
[53] Hu W, Ha YD, Bobaru F, Silling SA. The formulation and computation of the nonlocal j-integral in bond-based peridynamics. International journal of fracture 2012;176:195–206.
[54] Madenci E, Oterkus S. Ordinary state-based peridynamics for plastic deformation according to von mises yield criteria with isotropic hardening. Journal of the Mechanics and Physics of Solids 2016;86:192–219.
[55] Kakogiannou E, Sanavia L, Nicot F, Darve F, Schrefler BA. A porous media finite element approach for soil instability including the second-order work criterion. Acta Geotechnica 2016;11(4):805–25.
[56] Hill R. A general theory of uniqueness and stability in elastic-plastic solids. Journal of the Mechanics and Physics of Solids 1958;6(3):236–49.
[57] Gabriel E, Fagg GE, Bosilca G, Angskun T, Dongarra JJ, Squyres JM, et al. Open mpi: Goals, concept, and design of a next generation mpi implementation. In: Recent Advances in Parallel Virtual Machine and Message Passing Interface: 11th European PVM/MPI Users' Group Meeting Budapest, Hungary, September 19-22, 2004. Proceedings 11. Springer; 2004, p. 97–104.
[58] Zienkiewicz O, Paul D, Chan A. Unconditionally stable staggered solution procedure for soil-pore fluid interaction problems. International Journal for Numerical Methods in Engineering 1988;26(5):1039–55.
[59] Simoni L, Schrefler B. A staggered finite-element solution for water and gas flow in deforming porous media. Communications in applied numerical methods 1991;7(3):213–23.

[60] Armero F, Simo J. A new unconditionally stable fractional step method for non-linear coupled thermomechanical problems. International Journal for numerical methods in Engineering 1992;35(4):737–66.
[61] Turska E, Wisniewski K, Schrefler B. Error propagation of staggered solution procedures for transient problems. Computer methods in applied mechanics and engineering 1994;114(1-2):177–88.
[62] Schrefler B, Simoni L, Turska E. Standard staggered and staggered newton schemes in thermo-hydro-mechanical problems. Computer methods in applied mechanics and engineering 1997;144(1-2):93–109.
[63] Fredlund DG, Rahardjo H. Soil mechanics for unsaturated soils. John Wiley & Sons; 1993.
[64] Van Genuchten MT. A closed-form equation for predicting the hydraulic conductivity of unsaturated soils. Soil science society of America journal 1980;44(5):892–8.
[65] Song X, Borja RI. Finite deformation and fluid flow in unsaturated soils with random heterogeneity. Vadose Zone Journal 2014;13(5).
[66] Cao J, Jung J, Song X, Bate B. On the soil water characteristic curves of poorly graded granular materials in aqueous polymer solutions. Acta Geotechnica 2018;13(1):103–16.
[67] Niu WJ, Ye WM, Song X. Unsaturated permeability of gaomiaozi bentonite under partially free-swelling conditions. Acta Geotechnica 2020;15(5):1095–124.
[68] Ni T, Pesavento F, Zaccariotto M, Galvanetto U, Zhu QZ, Schrefler BA. Hybrid fem and peridynamic simulation of hydraulic fracture propagation in saturated porous media. Computer Methods in Applied Mechanics and Engineering 2020;366:113101.
[69] De Borst R. Simulation of strain localization: a reappraisal of the cosserat continuum. Engineering computations 1991;8(4):317–32.
[70] De Borst R, Mühlhaus HB. Gradient-dependent plasticity: formulation and algorithmic aspects. International journal for numerical methods in engineering 1992;35(3):521–39.
[71] Hughes TJ. The finite element method: linear static and dynamic finite element analysis. Courier Corporation; 2012.
[72] De Borst R, Crisfield MA, Remmers JJ, Verhoosel CV. Nonlinear finite element analysis of solids and structures. John Wiley & Sons; 2012.
[73] Belytschko T, Liu WK, Moran B, Elkhodary K. Nonlinear finite elements for continua and structures. John wiley & sons; 2014.
[74] Simo JC, Hughes TJ. Computational inelasticity; vol. 7. Springer Science & Business Media; 1998.
[75] Borja RI. Plasticity; vol. 2. Springer; 2013.
[76] Chen L, Fathi F, de Borst R. Hydraulic fracturing analysis in fluid-saturated porous medium. International Journal for Numerical and Analytical Methods in Geomechanics 2022;.